\documentclass{amsart}

\usepackage[latin1]{inputenc}
\usepackage{subfigure}
\usepackage{color}
\usepackage{amsmath, latexsym, amssymb,longtable}
\usepackage{graphicx}
\usepackage{amsfonts}
\usepackage{graphics}
\usepackage{latexsym}
\usepackage{amscd}
\usepackage{pict2e}
\usepackage[all]{xy}
\usepackage{bm}
\usepackage{enumerate}
\usepackage{todonotes}
\usepackage{mathrsfs}
\usepackage{shuffle}
\usepackage[mathscr]{euscript}

\newtheorem{theorem}{Theorem}[section]

\newtheorem{conjecture}[theorem]{Conjecture}
\newtheorem{definition}[theorem]{Definition}

\numberwithin{equation}{section}
\numberwithin{figure}{section}
\numberwithin{table}{section}

\DeclareMathOperator{\symgrp}{\mathfrak{S}}

\DeclareMathOperator{\leg}{\text{leg}}
\DeclareMathOperator{\Des}{Des}

\DeclareMathOperator{\maj}{maj}

\DeclareMathOperator{\inv}{inv}

\DeclareMathOperator{\SYT}{SYT}
\DeclareMathOperator{\NSym}{NSym}
\DeclareMathOperator{\FQSym}{FQSym}
\DeclareMathOperator{\Sym}{Sym}
\DeclareMathOperator{\QSym}{QSym}
\DeclareMathOperator{\SSRCT}{SSRCT}
\DeclareMathOperator{\SRCT}{SRCT}
\DeclareMathOperator{\South}{South}
\DeclareMathOperator{\SSYT}{SSYT}
\DeclareMathOperator{\Ind}{Ind}

\DeclareMathOperator{\asc}{asc}

\newlength\cellsize 
\savebox2{%
\begin{picture}(12,12)
\put(0,0){\line(1,0){12}}
\put(0,0){\line(0,1){12}}
\put(12,0){\line(0,1){12}}
\put(0,12){\line(1,0){12}}
\end{picture}}
\newcommand\cellify[1]{\def\thearg{#1}\def\nothing{}%
\ifx\thearg\nothing
\vrule width0pt height\cellsize depth0pt\else
\hbox to 0pt{\usebox2\hss}\fi%
\vbox to 12\unitlength{
\vss
\hbox to 12\unitlength{\hss$#1$\hss}
\vss}}
\newcommand\tableau[1]{\vtop{\let\\=\cr
\setlength\baselineskip{-16000pt}
\setlength\lineskiplimit{16000pt}
\setlength\lineskip{0pt}
\halign{&\cellify{##}\cr#1\crcr}}}
\savebox3{%
\begin{picture}(12,12)
\put(0,0){\line(1,0){12}}
\put(0,0){\line(0,1){12}}
\put(12,0){\line(0,1){12}}
\put(0,12){\line(1,0){12}}
\end{picture}}
\newcommand\expath[1]{%
\hbox to 0pt{\usebox3\hss}%
\vbox to 12\unitlength{
\vss
\hbox to 12\unitlength{\hss$#1$\hss}
\vss}}

\newcommand{\dI}{\mathcal{I}^*}
\newcommand{\I}{\mathcal{I}}

\author{ Sarah K. Mason}
\title{Recent trends in quasisymmetric functions}
\address{Department of Mathematics, Wake Forest University}

\begin{document}
\maketitle

 \paragraph{Abstract.}  
 
 This article serves as an introduction to several recent developments in the study of quasisymmetric functions.  The focus of this survey is on connections between quasisymmetric functions and the combinatorial Hopf algebra of noncommutative symmetric functions, appearances of quasisymmetric functions within the theory of Macdonald polynomials, and analogues of symmetric functions.  Topics include the significance of quasisymmetric functions in representation theory (such as representations of the $0$-Hecke algebra), recently discovered bases (including analogues of well-studied symmetric function bases), and applications to open problems in symmetric function theory.
 
 \tableofcontents
 
\section{Introduction}\label{sec:intro}

Quasisymmetric functions first appeared in the work of Stanley~\cite{Sta72} and were formally developed in Gessel's seminal article on multipartite $P$-partitions~\cite{Ges84}.  Since their introduction, their prominence in the field of algebraic combinatorics has continued to grow.  The number of recent developments in the study of quasisymmetric functions is far greater than would be reasonable to contain in this brief article; because of this, we choose to focus on a selection of subtopics within the theory of quasisymmetric functions.  This article is skewed toward bases for quasisymmetric functions which are closely connected to Macdonald polynomials and the combinatorial Hopf algebra of noncommutative symmetric functions.  A number of very interesting subtopics are therefore excluded from this article, including Stembridge's subalgebra of peak quasisymmetric functions~\cite{Ste97} and its associated structure (\cite{BMSvW00},~\cite{BMSvW02},~\cite{BHvW03},~\cite{BilLiu00},~\cite{Ehr96},~\cite{Li16}), Ehrenborg's flag quasisymmetric function of a partially ordered set~\cite{Ehr96}, colored quasisymmetric functions ~\cite{HsiKar11, HsiPet10}, and type B quasisymmetric functions~\cite{HsiPet06},~\cite{Pet05},~\cite{Pet07},~\cite{Cho01}.  This article also does not have the scope to address connections to probability theory such as riffle shuffles~\cite{Sta01}, random walks on quasisymmetric functions~\cite{HerHsi09}, or a number of other fascinating topics.  Hopefully this article will inspire the reader to learn more about quasisymmetric functions and explore these topics in greater depth.

This article begins with an overview of symmetric functions.  There are a number of excellent introductions to the subject including Fulton~\cite{Ful97}, Macdonald~\cite{Mac95}, Sagan~\cite{Sag01}, and Stanley~\cite{Sta99}.  The remainder of Section~\ref{sec:intro} deals with the genesis of quasisymmetric functions and several important bases.  Section~\ref{sec:rep} discusses the significance of quasisymmetric functions in algebra and representation theory, while Section~\ref{sec:mac} explores connections to Macdonald polynomials.  A number of recently introduced bases for quasisymmetric functions are described in Section~\ref{sec:bases}.  Section~\ref{sec:sym} is devoted to interactions with symmetric functions.

\subsection{Basic definitions and background on symmetric functions}{\label{sec:background}}

Recall that a \emph{permutation} of the set $[n] \colon=\{1,2, \hdots , n\}$ is a bijection from the set $\{1,2, \hdots , n\}$ to itself.  The group of all permutations of an $n$-element set is denoted $\mathfrak{S}_n$.  Let $\pi= \pi_1 \pi_2 \cdots \pi_n \in \mathfrak{S}_n$ denote a permutation written in one-line notation.  If $\pi_i > \pi_{i+1}$, then $i$ is a \emph{descent} of $\pi$.  If $\pi_i > \pi_j$ and $1 \le i < j \le n$, then the pair $(i,j)$ is an \emph{inversion} of $\pi$.  The \emph{sign} of a permutation $\pi$ (denoted $(-1)^{\pi}$) is the number of inversions of $\pi$.

Let $\mathbb{C}[x_1,x_2, \hdots , x_n]$ be the polynomial ring over the complex numbers $\mathbb{C}$ on a finite set of variables $\{x_1, x_2, \hdots, x_n \}$.  A permutation $ \pi \in \mathfrak{S}_n$ acts naturally on $f(x_1,x_2, \hdots , x_n) \in \mathbb{C}[x_1, x_2, \hdots , x_n]$ by $$\pi f(x_1, x_2, \hdots , x_n) = f(x_{\pi_1}, x_{\pi_2}, \hdots x_{\pi_n}).$$  

\begin{definition}
The ring of \emph{symmetric functions in $n$ variables} (often denoted by $\Lambda_n$ or $\Sym_n$) is the subring of $\mathbb{C}[x_1,x_2, \hdots , x_n]$ consisting of all polynomials invariant under the above action for all permutations in $\mathfrak{S}_n$.  
\end{definition}

This notion can be further extended to the ring $Sym$ of symmetric functions in infinitely many variables.  A \emph{symmetric function} $f \in \Sym$ is a formal power series $f \in \mathbb{C}[[X]]$ (with infinitely many variables $X=\{x_1,x_2, \hdots \}$) such that $f(x_1, x_2, \hdots )=f(x_{\pi_1}, x_{\pi_2}, \hdots)$ for every permutation of the positive integers. 

A \emph{partition} $\lambda = (\lambda_1, \lambda_2, \hdots , \lambda_{\ell})$ of a positive integer $n$ is a weakly decreasing sequence of positive integers which sum to $n$.  The elements $\lambda_i$ of the sequence are called the \emph{parts} and the number of parts is called the \emph{length} of the partition (denoted $\ell(\lambda)$).  We write $\lambda \vdash n$ (or $|\lambda| = n$) to denote ``$\lambda$ is a partition of $n$''.

\begin{figure}
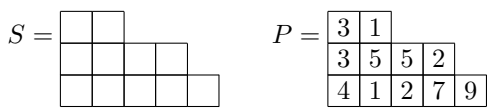

$$S = \tableau{ {} & {} \\ {} & {} & {} & {} \\ {} & {} & {} & {} & {} } \qquad P=\tableau{ 3 & 1 \\ 3 & 5 & 5 & 2 \\ 4 & 1 & 2 & 7 & 9 }$$
\caption{$S$ is the Ferrers diagram for $(5,4,2)$ and $P$ is a filling of $(5,4,2)$.}
\label{fig:ferrers}
\end{figure}

Each partition $\lambda =( \lambda_1, \lambda_2, \hdots , \lambda_{\ell})$ of $n$ can be visualized as a \emph{Ferrers diagram}, which consists of $n$ squares (typically called \emph{cells}) arranged into left-justified rows so that the $i^{th}$ row from the bottom contains $\lambda_i$ cells.  (Note that we are using French notation so that we think of the cells as indexed by their position in the coordinate plane.  This means the cell $(i,j)$ is the cell in the $i^{th}$ column from the left and the $j^{th}$ row from the bottom.  A Ferrers diagram in English notation places the rows so that the $i^{th}$ row from the top contains $\lambda_i$ cells, aligning with matrix indexing.)  An assignment of positive integer entries to each of the cells in the Ferrers diagram of shape $\lambda$ is called a \emph{filling}.  (See Figure~\ref{fig:ferrers}.)

A \emph{composition} $\alpha = (\alpha_1, \alpha_2, \hdots , \alpha_{\ell})$ of a positive integer $n$ is a sequence of positive integers which sum to $n$.  (It is sometimes necessary to allow $0$ to appear as a part; a \emph{weak composition} of a positive integer $n$ is a sequence of \emph{non-negative} integers which sum to $n$.)  The elements $\alpha_i$ of the sequence are called the \emph{parts} and the number of parts is called the \emph{length} of the composition.  Write $\alpha \models n$ to denote ``$\alpha$ is a composition of $n$''.  The \emph{reverse}, $\alpha^*$, of a composition $\alpha$ is obtained by reversing the order of the entries of $\alpha$ so that the last entry of $\alpha$ is the first entry of $\alpha^*$, the second to last entry of $\alpha$ is the second entry of $\alpha^*$, and so on.  For example, if $\alpha = (4,1,3,3,2)$, then $\alpha^* = (2,3,3,1,4)$.

\begin{figure}
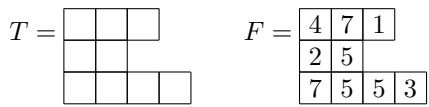

$$T=\tableau{ {} & {} & {} \\ {} & {} \\ {} & {} & {} & {} } \qquad F=\tableau{4 & 7 & 1 \\ 2 & 5 \\ 7 & 5 & 5 & 3}$$
\caption{$T$ is the composition diagram for $(4,2,3)$ and $F$ is a filling of $(4,2,3)$.}
\label{fig:Cdiagram}
\end{figure}

Each composition $\alpha = (\alpha_1, \alpha_2, \hdots , \alpha_{\ell})$ of $n$ can be visualized as a \emph{composition diagram}, which consists of $n$ cells arranged into left-justified rows so that the $i^{th}$ row from the bottom contains $\alpha_i$ cells, again using French notation.  (At times we will shift to English notation, but the reader may assume we are using French notation unless specified otherwise.)  An assignment of positive integer entries to each of the cells in the composition diagram of shape $\alpha$ is called a \emph{filling}.  See Figure~\ref{fig:Cdiagram} for an example of a composition diagram and a filling.

The \emph{refinement order} is a useful partial ordering on compositions.  We say $\alpha \prec \beta$ in the refinement ordering if $\beta$ can be obtained from $\alpha$ by summing adjacent parts of $\alpha$.  For example, $(2,4,1,1,3) \prec (2,5,4)$ and $(2,4,1,1,3) \prec (6,5)$ but $(2,5,4) \nprec (6,5)$ and $(6,5) \nprec (2,5,4)$.

Let $x^{\lambda}$ be the monomial $x_1^{\lambda_1} x_2^{\lambda_2} \cdots x_{\ell}^{\lambda_{\ell}}$.  For example, if $\lambda=(6,4,3,3,1)$, then $x^{\lambda} = x_1^6 x_2^4 x_3^3 x_4^3 x_5$.  One way to construct a symmetric function is to symmetrize such a monomial.  The \emph{monomial symmetric function} indexed by $\lambda$ is $$m_{\lambda} (X) = \sum x_{i_1}^{\lambda_1} x_{i_2}^{\lambda_2} \cdots x_{i_{\ell}}^{\lambda_{\ell}},$$ where the sum is over all distinct monomials with exponents $\lambda_1, \lambda_2, \hdots , \lambda_{\ell}$.

The monomial symmetric functions $\{ m_{\lambda} \mid \lambda \vdash n \}$ form a basis for the vector space $Sym^n$ of degree $n$ symmetric functions.  There are many elegant and useful bases for symmetric functions including three multiplicative bases obtained by describing the basis element $f_n$ and then setting $f_{\lambda} = f_{\lambda_1} f_{\lambda_2} \cdots f_{\lambda_{\ell}}$.  For example, the \emph{elementary symmetric functions} $\{ e_{\lambda} \mid \lambda \vdash n\}$ are defined by setting $e_n = m_{(1^n)}$, the \emph{complete homogeneous symmetric functions} $\{ h_{\lambda} \mid \lambda \vdash n\}$ are defined by setting $h_n = \displaystyle{\sum_{\lambda \vdash n} m_{\lambda}}$, and the \emph{power sum symmetric functions} $\{ p_{\lambda} \mid \lambda \vdash n\}$ are obtained by setting $p_n=m_{(n)}$.  Notice that the monomial symmetric functions are not multiplicative; that is, $m_{\lambda}$ is not necessarily equal to $m_{\lambda_1} m_{\lambda_2} \cdots m_{\lambda_{\ell}}$. 

Define a scalar product (a bilinear form $\langle f , g \rangle$ with values in $\mathbb{Q}$, sometimes referred to as an \emph{inner product}) on $\Sym$ by requiring that $$\langle h_{\lambda} , m_{\mu} \rangle = \delta_{\lambda \mu},$$ where $\delta$ is the Kronecker delta.  This means the complete homogeneous and monomial symmetric functions are dual to each other under this scalar product.  The power sums are orthogonal under this scalar product.  This means $\langle p_{\lambda}, p_{\mu} \rangle = \delta_{\lambda \mu} z_{\lambda}$, where $z_{\lambda} = 1^{m_1} (m_1!) 2^{m_2} (m_2)! \cdots k^{m_k} (m_k!)$ with $m_i$ equal to the number of times the value $i$ appears in $\lambda$.  For example, $$z_{(4,4,4,2,1,1)}=1^{2}2! 2^1 1! 3^0 0! 4^3 3! =1536.$$  (Note that $\frac{n!}{z_{\lambda}}$ counts the number of permutations of cycle type $\lambda$~\cite{Ber09, Sag01}.)  

Let $\omega \colon \Sym \rightarrow \Sym$ be the involution on symmetric functions defined by $\omega(e_n) = h_n$.  (Note that this implies $\omega(e_{\lambda}) = h_{\lambda}$ for all partitions $\lambda$.)  Then $\omega(p_n)=(-1)^{n-1} p_n$.

\subsection{Schur functions}

The \emph{Schur function} basis is one of the most important bases for symmetric functions due to its deep connections to representation theory and geometry as well as its combinatorial properties.  Schur functions are orthonormal under the scalar product described above, and can be defined in a number of different ways.  We begin with a combinatorial description, for which we will need several definitions.

\begin{figure}
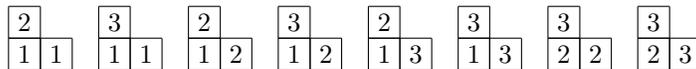

$$\tableau{2 \\ 1 & 1} \quad \tableau{3 \\ 1 & 1} \quad \tableau{2 \\ 1 & 2} \quad \tableau{3 \\ 1 & 2} \quad \tableau{2 \\ 1 & 3} \quad \tableau{3 \\ 1 & 3} \quad \tableau{3 \\ 2 & 2} \quad \tableau{3 \\ 2 & 3}$$
\caption{The set of all semi-standard Young tableaux of shape $(2,1)$ whose entries are in the set $\{1,2,3\}$}
\label{fig:ssyt}
\end{figure}

A filling of a partition diagram $\lambda$ in such a way that the row entries are weakly increasing from left to right and the column entries are strictly increasing from bottom to top is called a \emph{semi-standard Young tableau} of shape $\lambda$.  (See Figure~\ref{fig:ssyt}.)  The \emph{content} of such a filling is the composition $\alpha = (\alpha_1, \alpha_2, \hdots , \alpha_{k})$, where $\alpha_i$ is the number of times the entry $i$ appears in the filling.  The \emph{weight} of a semi-standard Young tableau $T$ is the monomial $x^T = x_1^{\alpha_1} x_2^{\alpha_2} \cdots x_k^{\alpha_k}$.  The set of all semistandard Young tableaux of shape $\lambda$ is denoted $\SSYT(\lambda)$.  A semi-standard Young tableau of shape $\lambda$ in which each entry from $1$ to $n$ (where $n = | \lambda|$) appears exactly once is called a \emph{standard Young tableau} and the set of all standard Young tableaux of shape $\lambda$ is denoted $\SYT(\lambda)$.

\begin{definition}
The \emph{Schur function} $s_{\lambda}(x_1,x_2, \hdots , x_n)$ is the generating function for the weights of all semi-standard Young tableaux of shape $\lambda$ with entries in the set $\{1,2, \hdots , n\}$; that is $$s_{\lambda}(x_1,x_2, \hdots , x_n) = \sum_{T \in \SSYT(\lambda)} x^T.$$  Here the sum is over all semi-standard Young tableaux whose entries are in the set $[n]$.  Extend this definition to infinitely many variables by allowing entries from the set of all positive integers.
\end{definition}

Figure~\ref{fig:ssyt} shows that $$s_{21}(x_1,x_2,x_3) = x_1^2x_2 + x_1^2x_3 + x_1 x_2^2 + 2 x_1 x_2 x_3 + x_1 x_3^2 + x_2^2x_3 + x_2 x_3^2.$$  Notice that $s_{(1^n)} = m_{(1^n)}=e_n$, since $s_{(1^n)}$ is constructed by filling a vertical column with positive integers so that no entries repeat, and $s_n=h_n$ since $s_n$ is constructed by filling a horizontal row with weakly increasing positive integers.  The Schur functions enjoy a large number of beautiful properties, including the celebrated Littlewood-Richardson formula for the coefficients appearing in the product of two Schur functions (which can also be computed algorithmically using the Remmel-Whitney Rule~\cite{RemWhi84}) and the Robinson-Schensted-Knuth Algorithm~\cite{Knu70,Sch61} which maps bijectively between matrices with finite non-negative integer support and pairs $(P,Q)$ of semi-standard Young tableaux of the same shape.  

The Schur functions were classically described as quotients involving the Vandermonde determinant and can be defined in a number of other ways.  One method of construction that can readily be generalized to other settings is through \emph{Bernstein creation operators}.

\begin{theorem}{\cite{Zel81}}
Define an operator ${{\bf B}}_m \colon \Sym^n \rightarrow \Sym^{m+n}$ by $${\bf B}_m :=\sum_{i \ge 0} (-1)^i h_{m+i}e_i^{\perp},$$ where $f^{\perp} \colon \Sym \rightarrow \Sym$ is an operator defined by $\langle fg,h \rangle = \langle g, f^{\perp}h \rangle$ for all $g,h \in Sym$.  Then for all tuples $\alpha \in \mathbb{Z}^m$, $$s_{\alpha} = {\bf B}_{\alpha_1} {\bf B}_{\alpha_2} \cdots {\bf B}_{\alpha_m} (1).$$
\end{theorem}

Note that this method for constructing Schur functions is more general than the combinatorial method described above because Bernstein creation operators define Schur functions indexed by tuples of non-negative integers (weak compositions) rather than just partitions.  

Schur functions appear in many areas of mathematics beyond combinatorics.  They correspond to characters of irreducible representations of the general linear group.  Their multiplicative structure describes the cohomology of the Grassmannian of subspaces of a vector space.  See the comprehensive texts by Fulton~\cite{Ful97} and Sagan~\cite{Sag01} for more details about Schur functions and their roles in combinatorics, representation theory, and geometry.

\subsection{Quasisymmetric functions}

The ring $\Sym$ of symmetric functions is contained inside a larger ring of \emph{quasisymmetric functions}, denoted by $\QSym$, which can be thought of as all bounded degree formal power series $f$ on an infinite alphabet $x_1, x_2, \hdots$ such that the coefficient of $x_{1}^{\alpha_1} x_{2}^{\alpha_2} \cdots x_{k}^{\alpha_k}$ in $f$ is equal to coefficient of $x_{j_1}^{\alpha_1} x_{j_2}^{\alpha_2} \cdots x_{j_k}^{\alpha_k}$ in $f$ for any sequence of positive integers $1 \le j_1 < j_2 < \cdots < j_k$ and any composition $(\alpha_1,\alpha_2, \hdots , \alpha_k)$.  It is often convenient to restrict to $n$ variables so that $f \in \QSym_n$ if and only if the coefficient of $x_{1}^{\alpha_1} x_{2}^{\alpha_2} \cdots x_{k}^{\alpha_k}$ in $f$ is equal to coefficient of $x_{j_1}^{\alpha_1} x_{j_2}^{\alpha_2} \cdots x_{j_k}^{\alpha_k}$ in $f$ for any sequence of positive integers $1 \le j_1 < j_2 < \cdots < j_k \le n$ and any composition $(\alpha_1,\alpha_2, \hdots , \alpha_k)$.

For example, the polynomial $$f(x_1,x_2,x_3) = x_1^3x_2^5+x_1^3 x_3^5 + x_2^3 x_3^5$$ is in $\QSym_3$, as is $$g(x_1,x_2,x_3) = x_1^3x_2^5+x_1^3 x_3^5 + x_2^3 x_3^5 + x_1^5x_2^3 + x_1^5 x_3^3 + x_2^5 x_3^3,$$ but $$h(x_1,x_2,x_3) = x_1^3 x_2^5 + x_1^3 x_3^5 + x_2^5 x_3^3$$ is not quasisymmetric in three variables since $x_2^3x_3^5$ does not appear, and neither do $x_1^5x_2^3$ and $x_1^5x_3^3$.  The quasisymmetric functions in $n$ variables are precisely the functions invariant under a quasi-symmetrizing action of the symmetric group $\mathfrak{S}_n$ introduced by Hivert~\cite{Hiv98,Hiv00}.

The origins of quasisymmetric functions first appeared in Stanley's work on $P$-partitions~\cite{Sta72}.  Gessel introduced the ring of quasisymmetric functions through his generating functions for Stanley's $P$-partitions~\cite{Ges84}.  (See Gessel~\cite{Ges16} for a historical survey of $P$-partitions.)  To be precise, let $[m]$ be the set $\{1,2, \hdots , m\}$ and let $X$ be an infinite totally ordered set such as the positive integers.  A \emph{partially ordered set} (or \emph{poset}), $(P, <_P)$, is a set of elements $P$ and a partial ordering $\le_P$ satisfying:
\begin{itemize}
\item reflexivity ($\forall x \in P, x \le_P x$), 
\item antisymmetry (if $x \le_P y$ and $y \le_P x$, then $x=y$), and 
\item transitivity (if $x \le_P y$ and $y \le_Pz$, then $x \le_P z$).  
\end{itemize}
Write $x <_P y$ if $x \le_P y$ and $x \not= y$.  A $P$-partition of a poset $P$ with elements $[m]$ and partial order $\le_P$ is a function $f \colon [m] \rightarrow X$ such that:
\begin{enumerate}
\item $i <_P j$ implies $f(i)$ is less than or equal to $f(j)$, and
\item $i <_P j$ and $i > j$ (under the usual ordering on integers) implies $f(i)$ is strictly less than $f(j)$.
\end{enumerate}


Each permutation $\pi$ corresponds to a totally ordered poset $P_{\pi}$ where $\pi_1 <_{\pi} \pi_2 <_{\pi} \cdots <_{\pi} \pi_m$.  It is these permutation posets that are used to construct Gessel's \emph{fundamental quasisymmetric functions $F_{\alpha}$}.  To do this, give each $P_{\pi}$-partition $f$ a \emph{weight} $x^f=\prod x_{f(i)}$ and sum the weights over all $P_{\pi}$-partitions.  

For example, let $\pi=312$ (written in one-line notation).  Condition (1) implies that $f(3) \le f(1) \le f(2)$.  Condition (2) implies $f(3)<f(1)$.   Therefore $f(3) < f(1) \le f(2)$, and the following table depicts the $P_{\pi}$-partitions  involving the subset $\{1,2,3 \}$ of $X$.

\begin{center}
\begin{tabular}{|c|c|c|c|}
\hline
f(3) & f(1) & f(2) & $x^f$ \\
\hline
1 & 2 & 2 & $x_1x_2^2$ \\
\hline
1 & 2 & 3 & $x_1x_2x_3$ \\
\hline
1 & 3 & 3 & $x_1x_3^2$ \\
\hline
2 & 3 & 3 & $x_2x_3^2$ \\
\hline
\end{tabular}
\end{center}

Therefore the fundamental quasisymmetric function corresponding to $\pi=(3,1,2)$ and restricted to three variables is $$ x_1x_2^2+x_1x_2x_3+x_1x_3^2+x_2x_3^2.$$  Note that this function depends only on the descent set of $\pi$ and not on the underlying permutation; therefore the indexing set for these generating functions is the set of all subsets of the set $[m-1]$ together with the number $m$ to indicate the degrees of the monomials.  We use an ordered pair consisting of a capital letter together with $m$ when indexing by sets.  In this paper we use Greek letters to denote compositions, but in other places, such as~\cite{GKLLRT95}, capital letters are used.

\begin{definition}{\cite{Ges84}}
Let $L$ be a subset of $[m-1]$.  Then
$$F_{L,m} = \sum_{\substack{{i_1 \le i_2 \le \cdots \le i_m}\\{i_j < i_{j+1}\; {\rm if} \; j \in L}}}x_{i_1} x_{i_2} \cdots x_{i_m}$$
\end{definition}

Every subset $L=\{L_1, L_2, \hdots , L_k\}$ of the set $[m-1]$ corresponds to a composition $\beta(L)=(L_1, L_2-L_1, \cdots, m-L_k)$ and therefore the fundamental quasisymmetric functions are often indexed by compositions rather than sets.    The fundamental quasisymmetric functions are homogeneous of degree $m$; the value $m$ is apparent when the index is a composition $\alpha$ since $|\alpha|=m$.  When the index is given by a subset this value $m$ must be specified.  For example, if $m=4$ then \begin{align*} F_{\{2,3\},4 }(x_1, x_2, x_3, x_4) &= \displaystyle{\sum_{i_1 \le i_2 < i_3 < i_4} x_{i_1} x_{i_2} x_{i_3} x_{i_4}} \\ &= x_1^2 x_2 x_3 +  x_1^2 x_2 x_4+ x_1^2 x_3 x_4+ x_1 x_2 x_3 x_4+ x_2^2 x_3 x_4.  \end{align*} whereas if $m=5$ then $$F_{\{2,3\},5 }(x_1, x_2, x_3, x_4) = \sum_{i_1 \le i_2 < i_3 < i_4 \le i_5} x_{i_1} x_{i_2} x_{i_3} x_{i_4} x_{i_5} \hspace*{1in}$$  $$\hspace*{.5in}= x_1^2x_2x_3^2 + x_1^2x_2 x_3 x_4 + x_1^2 x_2 x_4^2 + x_1^2 x_3 x_4^2 + x_1 x_2 x_3 x_4^2 + x_2^2 x_3 x_4^2.$$  

The \emph{monomial quasisymmetric function}, $M_{\alpha}$, in infinitely many variables $\{x_1, x_2, \hdots \}$ and indexed by the composition $\alpha = (\alpha_1, \alpha_2, \hdots , \alpha_k)$, is obtained by quasi-symmetrizing the monomial $x^{\alpha} =x_1^{\alpha_1} x_2^{\alpha_2} \cdots x_k^{\alpha_k}$.  That is, $$M_{\alpha} (x_1, x_2, \hdots ) = \sum_{i_1 < i_2 < \cdots < i_k} x_{i_1}^{\alpha_1} x_{i_2}^{\alpha_2} \cdots x_{i_k}^{\alpha_k}.$$  This definition can be restricted to finitely many variables by requiring that $i_k \le n$.  Note that $m_{\lambda} = \sum_{\tilde{\alpha}=\lambda} M_{\alpha},$ where $\tilde{\alpha}$ is the partition obtained by rearranging the parts of $\alpha$ into weakly decreasing order.  Every fundamental quasisymmetric function can be written as a positive sum of monomial quasisymmetric functions as follows: $$F_{\alpha} = \sum_{\beta \preceq \alpha} M_{\beta}.$$

The Schur functions decompose into a positive sum of fundamental quasisymmetric functions; to describe this decomposition we need one additional definition.  Each standard Young tableau $T$ has an associated \emph{descent set} $D(T)$ given by $i \in D(T)$ if and only if $i+1$ appears in a higher row of $T$ than $i$.  Then \begin{equation}{\label{Schur:fund}}
s_{\lambda} = \sum_{T \in \SYT(\lambda)} F_{D(T),|\lambda| }.
\end{equation}

For example, if $\lambda=(3,2)$, then the standard Young tableaux of shape $(3,2)$ are shown in Figure~\ref{fig:syta}, with respective descent sets $\{ 3 \}, \{ 2, 4 \}, \{ 2 \} , \{ 1 , 4\},$ and $\{1,3\}$.  Therefore \begin{align*}
s_{3,2} &= F_{ \{3\}, 5 } + F_{\{2,4 \}, 5 } + F_{\{2\},5} + F_{ \{1,4\},5} + F_{ \{1,3\},5} \\
&= F_{32}  +  F_{221} + F_{23} + F_{131} + F_{122}.
\end{align*}

\begin{figure}
$$\tableau{4 & 5 \\ 1 & 2 & 3} \qquad \tableau{3 & 5 \\ 1 & 2 & 4} \qquad \tableau{3 & 4 \\ 1 & 2 & 5} \qquad \tableau{2 & 5 \\ 1 & 3 & 4} \qquad \tableau{2 & 4 \\ 1 & 3 & 5}$$
\caption{The five standard Young tableau of shape $(3,2)$.}
{\label{fig:syta}}
\end{figure}

Valuable information can be gained about symmetric functions by examining their expansion into quasisymmetric functions, especially into the fundamental quasisymmetric functions.  For example, a symmetric function is said to be \emph{Schur positive} if it can be written as a positive sum of Schur functions.  Schur positivity is important because of its deep connection to representations of the symmetric group.  Assaf's recently developed paradigm called \emph{dual equivalence}~\cite{Ass15} provides machinery to prove that a function is Schur positive based on its expansion into the fundamentals and their connection to objects called dual equivalence graphs.  The \emph{Eulerian quasisymmetric functions}~\cite{ShaWac10} are defined as sums of the fundamental quasisymmetric functions indexed by certain permutation statistics.  Eulerian quasisymmetric functions are in fact always symmetric.  Their generating functions are deeply connected to Euler's exponential generating functions for the Eulerian polynomials.  Eulerian quasisymmetric functions can also be used to refine a number of classical results on permutation statistics.  We will not be able to address these topics in this brief survey article but encourage the interested reader to see~\cite{Ass15} and~\cite{ShaWac10} for details.

\section{Algebra and representation theory}{\label{sec:rep}}

Even before they were formally defined, quasisymmetric functions appeared naturally in algebraic settings.  The \emph{Leibniz-Hopf algebra} is the free associative algebra over the integers which in fact is isomorphic to the algebra of noncommutative symmetric functions, which we shall define in Section~\ref{sec:Hopf}.  In 1972, Ditters claimed~\cite[Proposition 2.2]{Dit72} that the Leibniz-Hopf algebra is a free commutative algebra over the integers.  This statement was later referred to as the \emph{Ditters Conjecture} due to an error in the original proof, and was then proved by Hazewinkel~\cite{Haz01, Haz10} using combinatorial techniques and later by Baker and Richter ~\cite{BakRic08} using methods from algebraic topology.  Malvenuto and Reutenauer~\cite[Corollary 2.2]{MalReu95} prove that $\QSym$ is a free module over $\Sym$.  In fact, Garsia and Wallach~\cite{GarWal03} further show that the quotient $\QSym_n$ over $\Sym_n$ has dimension $n!$.  Aval and Bergeron~\cite{ABB04} prove that the quotient of $\mathbb{Z}[x_1, x_2, \hdots , x_n]$ modulo quasisymmetric functions in $n$ variables with no constant term has Hilbert series $\sum C_n t^n$ where $C_n$ is the $n^{th}$ Catalan number (see also~\cite{AvaBer03} for the case with infinitely many variables).  An explicit basis for the quotient space is given in~\cite{LauMas10}.

The ring of quasisymmetric functions ($\QSym$) is an important example of a combinatorial Hopf algebra (discussed in Section~\ref{sec:Hopf}).  $\QSym$ is closely connected to Solomon's descent algebra (described in Section~\ref{sec:descent}) and plays a role in representations of the $0$-Hecke algebra (see Section~\ref{subsec:Hecke}).

\subsection{Combinatorial Hopf algebras}{\label{sec:Hopf}}

The following definitions, leading to the description of a combinatorial Hopf algebra, closely follow the expositions in~\cite{GriRei14} and~\cite{LMvW13}.

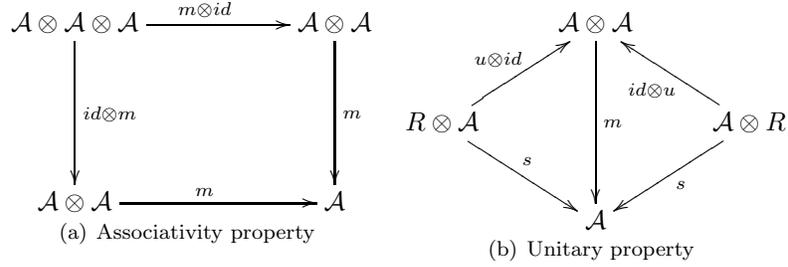
\begin{figure}{\label{fig:alg}}
\subfigure[Associativity property]{
\xymatrix{
\mathcal{A} \otimes \mathcal{A} \otimes \mathcal{A} \ar[rr]^{m \otimes id} \ar[dd]^{id \otimes m}  & & \mathcal{A} \otimes \mathcal{A} \ar[dd]^m  \\ \\
 \mathcal{A} \otimes \mathcal{A} \ar[rr]^m & & \mathcal{A}  \\ }
 }
 \subfigure[Unitary property]{
 \xymatrix{
 & \mathcal{A} \otimes \mathcal{A} \ar[dd]^m \\
 R \otimes \mathcal{A} \ar[ur]^{u \otimes id} \ar[rd]^{s} & &  \mathcal{A} \otimes R \ar[lu]^{id \otimes u} \ar[dl]^s \\
 & \mathcal{A}
 }
 }
 \caption{Commutative diagrams for algebras}
 \end{figure}

Let $R$ be a commutative ring with an identity element.  An \emph{associative algebra} over $R$ is an $R$-module $\mathcal{A}$ together with a \emph{product} (or \emph{multiplication}) $m \colon \mathcal{A} \otimes \mathcal{A} \rightarrow \mathcal{A}$ and a \emph{unit} $u \colon R \rightarrow \mathcal{A}$ satisfying associativity ($m(m(a,b),c) = m(a, m(b,c))$) and a unitary property which implies that the unit map commutes with scalar multiplication.  To be precise, $m$ and $u$ are $R$-linear maps such that if $id$ is the identity map on $\mathcal{A}$ and $s$ is scalar multiplication, then the diagrams in Figure~\ref{fig:alg} commute.

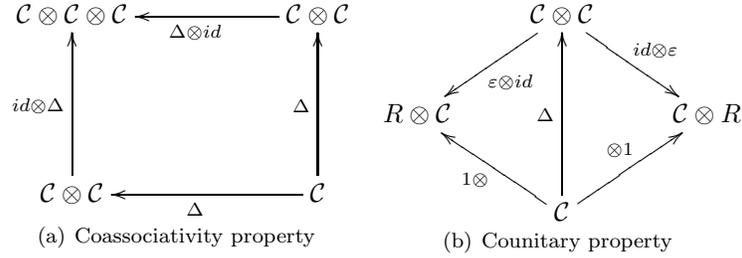
\begin{figure}
\subfigure[Coassociativity property]{
\xymatrix{
\mathcal{C} \otimes \mathcal{C} \otimes \mathcal{C}  & & \mathcal{C} \otimes \mathcal{C} \ar[ll]^{\Delta \otimes id} \\ \\
 \mathcal{C} \otimes \mathcal{C} \ar[uu]^{id \otimes \Delta} & & \mathcal{C} \ar[ll]^{\Delta} \ar[uu]^{\Delta} \\ }
 }
 \subfigure[Counitary property]{
 \xymatrix{
 & \mathcal{C} \otimes \mathcal{C} \ar[dr]^{id \otimes \varepsilon} \ar[dl]^{\varepsilon \otimes id} \\
 R \otimes \mathcal{C}  & &  \mathcal{C} \otimes R \\
 & \mathcal{C} \ar[uu]^{\Delta} \ar[ur]^{\otimes 1} \ar[ul]^{1 \otimes}
 }
 }
 \caption{Commutative diagrams for coalgebras  (Here $\otimes 1$ applied to an element $c \in \mathcal{C}$ means $c \otimes 1$ and $1 \otimes $ applied to an element $c \in \mathcal{C}$ means $1 \otimes c$.)}{\label{fig:coalg}}
 \end{figure}

A \emph{coalgebra} over $R$ is an $R$-module $\mathcal{C}$ together with a \emph{coproduct} (or \emph{comultiplication}) $\Delta \colon \mathcal{C} \rightarrow \mathcal{C} \otimes \mathcal{C}$ and a \emph{counit} (or \emph{augmentation}) $\varepsilon \colon \mathcal{C} \rightarrow R$ satisfying a \emph{coassociativity} property and a \emph{counitary} property so that when the directions of the arrows in Figure~\ref{fig:alg} are reversed, the resulting diagrams (shown in Figure~\ref{fig:coalg}) commute.

An \emph{algebra morphism} is a map $f \colon \mathcal{A} \rightarrow \mathcal{A}'$ from an $R$-algebra $(\mathcal{A}, m , u)$ to another $R$-algebra $(\mathcal{A}', m', u')$ such that $$f \circ m = m' \circ (f \otimes f) \; {\rm and} \; f \circ u = u'.$$ A \emph{bialgebra} is an algebra $(\mathcal{B}, m ,u)$ and coalgebra $(\mathcal{B}, \Delta, \varepsilon)$ such that $\Delta$ and $\varepsilon$ are algebra homomorphisms.  A bialgebra $\mathcal{B}$ with coproduct $\Delta$ is said to be \emph{graded} if it decomposes into submodules $\mathcal{B}_0, \mathcal{B}_1, \hdots$ such that

\begin{enumerate}
\item $\mathcal{B} = \displaystyle{\oplus_{n \ge 0} \mathcal{B}_n}$,
\item $\mathcal{B}_i \mathcal{B}_j \subseteq \mathcal{B}_{i+j}$, and
\item $\Delta(\mathcal{B}_n) \subseteq \oplus_{i+j=n} \mathcal{B}_i \otimes \mathcal{B}_j$.
\end{enumerate}

\begin{definition}
A bialgebra $(\mathcal{H}, m ,u, \Delta, \varepsilon)$ is a \emph{Hopf algebra} if there exists a linear map $S \colon \mathcal{H} \rightarrow \mathcal{H}$ (called the \emph{antipode}) such that $$m \circ (S \otimes id) \circ \Delta = u \circ \varepsilon = m \circ (id \otimes S) \circ \Delta.$$
\end{definition}

A Hopf algebra $\mathcal{H}$ is said to be \emph{connected} if $\mathcal{H}_0=R$.  When the ground ring $R$ is in fact a field $K$, a \emph{character} (sometimes called a \emph{multiplicative linear functional}) of the Hopf algebra $\mathcal{H}$ is an algebra homomorphism from $\mathcal{H}$ to the field $K$.  A \emph{combinatorial Hopf algebra} is a graded connected Hopf algebra equipped with a character.

Gessel~\cite{Ges84} describes an \emph{internal (or inner) coproduct} which takes $\QSym^n$ (quasisymmetric functions of degree $n$) to $\QSym^n \otimes \QSym^n$.  This internal coproduct corresponds to the internal coproduct on symmetric functions, taking $p_n$ to $p_n \otimes p_n$.  Malvenuto and Reutenauer~\cite{MalReu95} introduce an \emph{outer coproduct} on $\QSym$ defined on the monomial quasisymmetric functions by $$\Delta (M_{(\beta_1, \beta_2, \hdots , \beta_k)}) = \sum_{i=0}^k M_{(\beta_1, \hdots , \beta_i)} \otimes M_{(\beta_{i+1}, \hdots , \beta_k)}.$$  For example, $$\Delta(M_{312}) = 1 \otimes M_{312} + M_3 \otimes M_{12} + M_{31} \otimes M_2 + M_{312} \otimes 1.$$  Restricting this coproduct to symmetric functions takes $p_n$ to $p_n \otimes 1 + 1 \otimes p_n$, and therefore this coproduct is different from Gessel's internal coproduct.  Malvenuto and Reutenauer~\cite{MalReu95} and Ehrenborg~\cite{Ehr96} independently discovered the antipode map on $\QSym$ (with respect to the outer coproduct), proving that $\QSym$ is a Hopf algebra.  See~\cite{Haz03} for a thorough introduction to the Hopf algebra structure of quasisymmetric functions.

If $V$ is an $R$-module, let $V^{\star}:=Hom(V,R)$ be its dual $R$-module.  Under certain finiteness conditions (which are true for the situations explored in this article), the duals of Hopf algebras are themselves Hopf algebras.  In the language of Hopf algebras and their duality, the Ditters conjecture states that the Leibniz-Hopf algebra is dual (as a Hopf algebra over the integers) to a free commutative algebra over the integers.

The dual to $\QSym$ is the ring (or algebra) of noncommutative symmetric functions, denoted $\NSym$.  We take a moment to briefly describe some of the structure of $\NSym$.  For a thorough introduction to the topic through the lens of quasi-determinants, please see~\cite{GKLLRT95}; we typically follow their notation conventions.  

$\NSym$ can be thought of as a free associative algebra $K \langle {\bf \Lambda}_1, {\bf \Lambda}_2, \hdots \rangle$ generated by an infinite sequence of noncommuting indeterminates $({\bf \Lambda}_k)_{k \ge 1}$ over a fixed field $K$ of characteristic $0$.  (We usually take $K$ to be $\mathbb{C}$, the complex numbers.)  The \emph{noncommutative elementary symmetric functions} are the indeterminates ${\bf \Lambda}_k$, and their generating function is $$\lambda(t) = \sum_{k \ge 0} t^k {\bf \Lambda}_k,$$ while the \emph{noncommutative complete homogeneous symmetric functions} $S_k$ are defined by their generating function \begin{equation}{\label{eqn:genfun}} \sigma(t) = \sum_{k \ge 0} t^k {\bf S}_k = \lambda(-t)^{-1}.\end{equation}  Note that this mirrors the relationship in $\Sym$ between the elementary and complete homogeneous symmetric functions, where if $H(t)=\sum_{n \ge 0} h_nt^n$ and $E(t) = \sum_{n \ge 0} e_nt^n$, then $H(t)E(-t)=1$.  Both of these basis analogues are multiplicative, meaning ${\bf S}_{\alpha} = {\bf S}_{\alpha_1} {\bf S}_{\alpha_2} \cdots {\bf S}_{\alpha_k}$ and ${\bf \Lambda}_{\alpha} = {\bf \Lambda}_{\alpha_1} {\bf \Lambda}_{\alpha_2} \cdots {\bf \Lambda}_{\alpha_k}$ for $\alpha = (\alpha_1, \alpha_2, \hdots \alpha_k)$.  The \emph{Ribbon Schur functions}, which form a basis dual to Gessel's fundamental basis for $\QSym$, can be defined through quasi-determinants.  Two different candidates for the noncommutative analogue of the power sum symmetric functions will be described in Section~\ref{sec:powersum}.  We use boldface letters for bases of $\NSym$, lowercase letters for bases of $\Sym$, and uppercase letters for bases of $\QSym$.

The \emph{forgetful map}, frequently denoted by $\chi \colon \NSym \rightarrow \Sym$, sends the basis element ${\bf S}_{\alpha}$ to the complete homogeneous symmetric function $h_{\alpha_1} h_{\alpha_2} \cdots h_{\alpha_{\ell(\alpha)}}$.  Essentially, the forgetful map ``forgets'' that the functions don't commute.  This map can then be extended linearly to all of $\NSym$ and is in fact a surjection (but clearly not a bijection) onto $\Sym$.

Let $C$ be a category of objects.  An object $T$ is a \emph{terminal object} for the category $C$ if for all objects $X \in C$ there exists a unique morphism $X \rightarrow T$.  Not every category has a terminal object, but if such a terminal object exists it is necessarily unique.  Aguilar, Bergeron, and Sottile~\cite{ABS06} introduce a canonical character $\zeta_Q$ on $\QSym$ and describe what it does to the monomial and fundamental quasisymmetric functions.  Equipped with this character, quasisymmetric functions are the terminal object in the category of combinatorial Hopf algebras. 

\begin{theorem}{\cite{ABS06}}{\label{ABSthm}}
If $\mathcal{H}$ is a combinatorial Hopf algebra with a character $\zeta$, then there exists a unique homomorphism from $(\mathcal{H}, \zeta)$ to $(\QSym, \zeta_Q)$ such that the homomorphism on characters induced by the Hopf algebra homomorphism sends $\zeta$ to $\zeta_Q$.
\end{theorem}

Theorem~\ref{ABSthm} helps to explain why quasisymmetric functions appear in so many different contexts throughout algebraic combinatorics.  Examples for which the connection is well-understood include Rota's Hopf algebra of isomorphism classes of finite graded posets~\cite{JonRot79} and the chromatic Hopf algebra of isomorphism classes of finite unoriented graphs~\cite{Sch94}.  Note that this mirrors the similar result stating that $\Sym$ is the terminal object in the category of cocommutative combinatorial Hopf algebras~\cite{ABS06}.  

\subsection{Solomon's descent algebra}{\label{sec:descent}}

Let $\mathbb{Z}\!\symgrp_n$ be the group ring of permutations $\symgrp_n$ over the integers and let $\displaystyle{\mathbb{Z}\!\symgrp = \oplus_{n \ge 0} \mathbb{Z}\!\symgrp_n}$ be the direct sum of $\mathbb{Z}\!\symgrp_n$ over all positive integers $n$.  For each $\sigma \in \symgrp_n$ let $\Des(\sigma)$ be the \emph{descent set} of $\sigma$ defined by $\Des(\sigma) = \{ i \mid 1 \le i \le n-1, \sigma(i)> \sigma(i+1) \}$.  To each subset $L$ of $\{1,2, \hdots , n-1\}$, associate an element $D_L$ of $\mathbb{Z}\symgrp_n$ as follows:  $$D_L = \sum_{\Des(\sigma) = L} \sigma.$$  The composition $\beta(L)=(L_1, L_2-L_1, \hdots , n-L_k)$ (where $k=|L|$) is frequently used to index this \emph{descent basis}.  For example, let $n=4$ and $L=\{2\}$.  Then $$D_L = D_{(2,2)}= 1324 + 1423 + 2314 + 2413 + 3412.$$

Let $\Sigma_n$ be the linear span of the elements $D_L$ and endow $\displaystyle{\Sigma:=\oplus_{n \ge 0} \Sigma_n \subseteq \mathbb{Z}\!\symgrp}$ with a ring structure by setting $\sigma \pi = 0$ if $\sigma \in \symgrp_n$ and $\pi \in \symgrp_m$ such that $m \not= n$.  $\Sigma$ is called \emph{Solomon's descent algebra}.  Solomon~\cite{Sol76} proves that $\Sigma$ is a subalgebra of $\mathbb{Z}\!\symgrp$.  Gessel~\cite{Ges84} shows that the algebra dual to the coalgebra $\QSym_n$ (endowed with Gessel's internal coproduct) is isomorphic to Solomon's descent algebra $\Sigma_n$.  




The set $\Sigma$ (just as a set, not as the descent algebra) also admits a Hopf algebra structure.  That is, Malvenuto and Reutenauer~\cite{Mal94,MalReu95} define a product and coproduct on $\mathbb{Z}\!\symgrp$ to prove that $\mathbb{Z}\!\symgrp$ is a Hopf algebra (called the \emph{Malvenuto-Reutenauer Algebra}) with $\Sigma$ as a Hopf subalgebra.  (This algebra is in fact isomorphic to the algebra $\FQSym$ of \emph{free quasi-symmetric functions}; see~\cite{DHT02} for details.)  Malvenuto and Reutenauer show (Theorem 3.3 in \cite{MalReu95}) that $\Sigma$ is dual to $\QSym$, with the descent basis $\{ D_{\alpha} \}$ of $\Sigma$ dual to the basis $\{ F_{\alpha} \}$.  This means that the product on one of these bases determines the coproduct on the other, and vice versa.  That is, if $$F_{\alpha} F_{\beta} = \sum_{\gamma \models |\alpha| + | \beta|} c^{\gamma}_{\alpha, \beta} F_{\gamma},$$ then comultiplication in $\Sigma$ is defined by $$\Delta_{\Sigma} (D_{\gamma}) = \sum c^{\gamma}_{\alpha, \beta} D_{\alpha} \otimes D_{\beta}.$$

This duality pairing also implies that the descent basis is isomorphic to the ribbon Schur basis for noncommutative symmetric functions since the ribbon Schur basis for $\NSym$ is dual to the fundamental basis for $\QSym$.  

\subsection{Representations of the $0$-Hecke algebra}{\label{subsec:Hecke}}

The representation theoretic significance of the fundamental quasisymmetric functions mirrors that of the Schur functions.  We first describe the symmetric function case for ease of comparison.  Recall that the symmetric group $\symgrp_n$ can be generated by adjacent transpositions $s_i=(i,i+1)$ for $1 \le i \le n-1$ satisfying the following relations:
$$s_i^2 = 1 \; {\rm for} \; 1 \le i \le n-1,$$
$$s_i s_{i+1} s_i = s_{i+1} s_i s_{i+1} \; {\rm for} \; 1 \le i \le n-2, \; \text{and}$$
$$s_is_j=s_j s_i  \; {\rm for} \;|i-j| \ge 2.$$

The Frobenius \emph{characteristic map} is a map from characters of the symmetric group $\mathfrak{S}_n$ to symmetric functions which are homogeneous of degree $n$.  The Schur functions are the images of irreducible characters.

The $0$-Hecke algebra is a $\mathbb{C}$-algebra generated by elements satisfying relations similar to the relations on the symmetric group generators described above.  That is, $H_n(0)$ is generated by elements $T_1, T_2, \hdots , T_{n-1}$ satisfying: $$T_i^2 = T_i \; {\rm for} \; 1 \le i \le n-1,$$
$$T_i T_{i+1} T_i = T_{i+1} T_i T_{i+1} \; {\rm for} \; 1 \le i \le n-2, \; \text{and}$$
$$T_i T_j = T_j T_i \; {\rm if} \; | i-j| \ge 2.$$

If $\sigma$ is a permutation in $\mathfrak{S}_n$ with reduced word $\sigma = s_{i_1} s_{i_2} \cdots s_{i_{\ell}}$, then define $T_{\sigma} \in H_n(0)$ by $$T_{\sigma} = T_{i_1} T_{i_2} \cdots T_{i_{\ell}}.$$  The $0$-Hecke algebra $H_n(0)$ is a specialization of the Hecke algebra $H_n(q)$ at $q=0$.  (See M\'{e}liot~\cite{Mel17} for further details on the Hecke algebra $H_n(q)$ and its relationship to the $0$-Hecke algebra $H_n(0)$.)  

Norton~\cite{Nor79} investigates the representation theory of $H_n(0)$ and proves that there are $2^{n-1}$ distinct irreducible representations of $H_n(0)$, indexed by compositions of $n$.  Let $\mathcal{G}_0 (H_n(0))$ be the Grothendieck group of finitely generated $H_n(0)$-modules and $\displaystyle{\mathcal{G} = \oplus_{n \ge 0} \mathcal{G}_0(H_n(0))}$ be the associated Grothendieck ring.  (See Carter~\cite{Car86} for a thorough account of the representation theory of the $0$-Hecke algebra and see Huang~\cite{Hua14,Hua15,Hua16} for recent connections with flag varieties, the Stanley-Reisner ring, and tableaux.)  Krob and Thibon~\cite{KroThi97} prove that $\mathcal{G}$ is isomorphic to the ring of quasisymmetric functions via a characteristic map $\mathcal{F} \colon \mathcal{G} \rightarrow \QSym$ called the \emph{quasisymmetric characteristic}.  Let $L_{\alpha}$ denote the irreducible representation of $H_n(0)$ corresponding to $\alpha$.  Then $\mathcal{F}$ sends $L_{\alpha}$ to the fundamental quasisymmetric function $F_{\alpha}$.

\begin{theorem}{\cite{DKLT96}}
The map $\mathcal{F}$, defined by $\mathcal{F}(L_{\alpha})=F_{\alpha}$, is a ring isomorphism between the Grothendieck group $\mathcal{G}$ of finite-dimensional representations of $H_n(0)$ and the ring of quasisymmetric functions.
\end{theorem}

The fundamental quasisymmetric functions therefore correspond to characters of irreducible representations of the $0$-Hecke algebra.  See Hivert~\cite{Hiv98, Hiv00} to view this through the lens of divided difference operators.  Similar representation theoretic interpretations can be ascribed to various other bases for quasisymmetric functions and will be discussed in the relevant sections.

\section{Macdonald polynomials}{\label{sec:mac}}

The Schur functions are uniquely determined by the following two requirements described on p.305 of Macdonald~\cite{Mac98}.

\begin{enumerate}
\item Let $\lambda$ be a partition.  Then $$s_{\lambda} = m_{\lambda} + \sum_{\mu < \lambda} K_{\lambda \mu} m_{\mu},$$ where $\mu \le \lambda$ if and only if $\mu_1 + \mu_2 + \cdots + \mu_j \le \lambda_1 + \lambda_2 + \cdots \lambda_j$ for all $j$.  (This partial ordering is called the \emph{dominance ordering}.)  Here the coefficients $K_{\lambda \mu}$ are called the \emph{Kostka numbers}, or \emph{Kostka coefficients}.
\item $\langle s_{\lambda}, s_{\mu} \rangle = \delta_{\lambda \mu}$.
\end{enumerate}

Macdonald~\cite{Mac88} generalized this construction to a two-parameter family of functions $P_{\lambda} = P_{\lambda}(q,t)$ in $\mathbb{Q}(q,t)$ characterized by the following two requirements.

\begin{enumerate}
\item Let $\lambda$ be a partition.  Then $P_{\lambda} = m_{\lambda} + \mbox{lower terms in dominance order}$.
\item $\langle P_{\lambda}, P_{\mu} \rangle_{q,t} = 0 \; {\rm if} \; \lambda \not= \mu,$ where $$\langle p_{\lambda} , p_{\mu} \rangle_{q,t} = \delta_{\lambda \mu} z_{\lambda} \prod_{i=1}^{\ell(\lambda)} \frac{1-q^{\lambda_i}}{1-t^{\lambda_i}}.$$
\end{enumerate}

Notice that when $q=t$, the scalar product reduces to $\langle p_{\lambda}, p_{\mu} \rangle_{q,q} = \delta_{\lambda \mu} z_{\lambda}$ and so $P_{\lambda} (q,q) = s_{\lambda}$.  Similarly, $P_{\lambda}(q,1)=m_{\lambda}$ and $P_{\lambda}(1,t) = e_{\lambda'}$.  (See p.324 of Macdonald~\cite{Mac98}.)  Therefore the Macdonald polynomials are simultaneous generalizations of several different symmetric function bases.  Macdonald polynomials also appear in connection with the Hilbert scheme of $n$ points in the plane~\cite{Hai01}.

There are several variations on the original definition of Macdonald polynomials, including the \emph{modified Macdonald polynomials}, $\tilde{H}_{\mu}$, obtained from $P_{\mu}$ by certain substitutions and motivated by their connection to the coefficients appearing in the Schur function expansion of Macdonald polynomials.  Haglund conjectured and Haglund, Haiman, and Loehr proved~\cite{Hag04a, Hag05, HHL04, HHL05} a combinatorial formula for the Macdonald polynomials $\tilde{H}_{\mu}$ using statistics on fillings of partition diagrams.  To describe this formula, we introduce several pertinent definitions.

Recall that a filling $\sigma \colon \mu \rightarrow \mathbb{Z}^+$ is a function from the cells of the diagram of a partition $\mu$ to the positive integers.  The \emph{reading word} of the filling is the word obtained by reading the entries of a filling from top to bottom, left to right.  

\begin{figure}
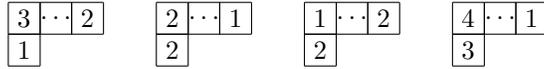

$$\tableau{3 &  \hdots & 2 \\ 1} \qquad \tableau{ 2 & \hdots & 1 \\ 2} \qquad \tableau{1 & \hdots & 2 \\ 2} \qquad \tableau{4 & \hdots & 1 \\ 3}$$
\caption{The first two triples are inversion triples; the third and fourth are not.}
\label{fig:invtrips}
\end{figure}

The major index and inversion statistic on permutations can be generalized to statistics on fillings of Ferrers diagrams.  Let $s$ be a cell in the partition diagram $\mu$ and let $\South(s)$ be the cell immediately below $s$ in the same column as $s$.  Define $$\Des(\sigma, \mu) = \{ s \in \mu \mid \sigma(s) > \sigma(\South(s))\}.$$  (No cell in the bottom row of $\mu$ can be in $\Des(\sigma, \mu)$.)  Let $\leg(s)$ be the number of cells above $s$ in the same column as $s$ and let $arm(s)$ be the number of cells to the right of $s$ in the same row as $s$.  Then $$\maj(\sigma,\mu) = \sum_{s \in \Des(\sigma,\mu)} (\leg(s) +1).$$  Let $u,v,w$ be three cells in the diagram of $\mu$ such that $u$ and $v$ are in the same row of $\mu$ with $v$ strictly to the right of $u$ and $w=\South(u)$ as shown: $$\tableau{u & \hdots & v \\ w}.$$  Any collection of three cells arranged in this way is called a \emph{triple}.  Define an orientation on the cells in a triple of a filling $\sigma$ of $\mu$ by starting with the cell containing the smallest entry and moving in a circular motion from smallest to largest.  (If two entries are equal, the one which appears first in the reading word is considered smaller.)  If the resulting orientation is counterclockwise, the triple is called an \emph{inversion triple}.  (See Figure~\ref{fig:invtrips}.)  Two cells $u,v$ in the bottom row are also considered an inversion triple if $v$ is strictly to the right of $u$ and $\sigma(u) > \sigma(v)$.  The total number of inversion triples in a filling $\sigma$ of a partition $\mu$ is denoted $\inv(\sigma, \mu)$.  

\begin{figure}
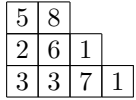

$$\tableau{5 & 8 \\ 2 & 6 & 1 \\ 3 & 3 & 7 & 1}$$
\caption{A filling of the partition $(4,3,2)$ with reading word $582613371$}
\label{fig:macfilling}
\end{figure}

For example, the filling in Figure~\ref{fig:macfilling} has descent set $\Des=\{ (1,3), (2,3), (2,2) \}$, where cells are indexed by (column, row) to mimic the $(x,y)$ Cartesian coordinates.  The major index for this filling is $1+1+2=4$ and $\inv(\sigma, \mu)=3+2=5$, since there are three inversion triples in the bottom row and two additional inversion triples.

\begin{theorem}{\cite{HHL04, HHL05}}{\label{thm:maccomb}}
Let $\mu$ be a partition of $n$.  Then $$\tilde{H}_{\mu}(X;q,t) = \sum_{\sigma \colon \mu \rightarrow \mathbb{Z}^+} x^{\sigma}  q^{\inv(\sigma, \mu)} t^{\maj(\sigma, \mu)}.$$
\end{theorem}

Theorem~\ref{thm:maccomb} provides a straightforward method for computing Macdonald polynomials.  This formula could potentially be used to find a product rule for Macdonald polynomials utilizing tableau constructions, although Yip recently found an elegant combinatorial rule for multiplying Macdonald polynomials~\cite{Yip12} using the alcove walk model introduced by Ram and Yip~\cite{RamYip11}.

\subsection{Quasisymmetric decomposition of Macdonald polynomials}

Macdonald polynomials can also be described as sums of fundamental quasisymmetric functions with coefficients in $q$ and $t$.

\begin{theorem}{\cite{HHL04,Hag08}}{\label{thm:HHL}}
Let $\mu$ be a partition of $n$.  Then 
$$\tilde{H}_{\mu}(X; q,t) = \sum_{\beta \in \mathfrak{S}_n} q^{\mathrm{\inv}(\beta,\mu)} t^{\mathrm{\maj}(\beta,\mu)}F_{\Des(\beta^{-1})},$$
where each permutation $\beta$ in the sum corresponds to the standard filling of $\mu$ with reading word $\beta$ and $\Des(\beta^{-1})$ is the usual descent set of the permutation $\beta^{-1}$ obtained by taking the inverse of $\beta$.
\end{theorem}

For example, if $\mu=(2,1)$, Table~\ref{table:Hmu} demonstrates that $$\tilde{H}_{21} (X;q,t) = F_3 + (q+t)F_{21} + (q+t)F_{12} + qtF_{111}.$$

\begin{table}{\label{table:Hmu}}
\begin{center}
\begin{tabular}{| c | c | c | c | c | c | c |}
\hline
Permutation $\beta$ & 123 & 132 & 213 & 231 & 312 & 321 \\
(Reading word of filling) & & & & & & \\
\hline
Filling of $\mu$ & $\tableau{1 \\ 2 & 3 }$ & $\tableau{1 \\ 3 & 2}$ & $\tableau{2 \\ 1 & 3}$ & $\tableau{2 \\ 3 & 1}$ & $\tableau{3 \\ 1 & 2}$ & $\tableau{3 \\ 2 & 1}$ \\
\hline
$\inv(\beta, \mu)$ & 0 & 1 & 0 & 1 & 0 & 1 \\
\hline
$\maj(\beta, \mu)$ & 0 & 0 & 1 & 0 & 1 & 1 \\
\hline
$\beta^{-1}$ & 123 & 132 & 213 & 312 & 231 & 321 \\
\hline
$\Des(\beta^{-1})$ & $\emptyset$ & 2 & 1 & 1 & 2 & 1,2 \\
\hline
\end{tabular}
\end{center}
\end{table}

This expansion of the Macdonald polynomials into fundamental quasisymmetric functions paves the way for new approaches to long-standing open questions.  For example, Macdonald~\cite{Mac88} conjectured that the coefficients in the expansion of $\tilde{H}_{\mu}$ into Schur functions are polynomials in $q$ and $t$ with nonnegative integer coefficients.  Haiman~\cite{Hai01} proved this by showing that $\tilde{H}_{\mu}$ is the bigraded Frobenius character of a doubly-graded $\mathfrak{S}_n$-module, but this approach did not provide an explicit combinatorial formula for the coefficients.  Assaf's dual equivalence~\cite{Ass15,Ass07} provides another potential approach to Schur positivity which makes use of the decomposition of a symmetric function into fundamental quasisymmetric functions.

The Hall-Littlewood polynomials are a one-parameter specialization of Macdonald polynomials introduced by Littlewood as a symmetric function realization of the Hall algebra~\cite{Lit61}.  Several different candidates for quasisymmetric Hall-Littlewood polynomials have recently been proposed.  See Hivert~\cite{Hiv00} for an analogue in $\NSym$ and its $\QSym$ companion, and see Novelli, Thibon, and Williams~\cite{NTW10} for a different noncommutative analogue.  Connections between these two approaches are studied in Novelli, Tevlin, and Thibon~\cite{NTT13}.  See also Haglund, Luoto, Mason, and van Willigenburg~\cite{HLMvW09} for another quasisymmetric analogue.

\subsection{Quasisymmetric Schur functions}

Haglund's formula (Theorem~\ref{thm:maccomb}) to generate the Macdonald polynomials using statistics on fillings of partition diagrams is generalized in~\cite{HHL08} to fillings of weak composition diagrams in order to generate the nonsymmetric Macdonald polynomials introduced and initially developed by Cherednik~\cite{Che95}, Macdonald~\cite{Mac95}, Opdam~\cite{Opd95}, and Sahi~\cite{Sah99}.  When these polynomials are specialized to $q=t=\infty$, the resulting polynomials, called \emph{Demazure atoms} due to their connections to Demazure characters, form a basis for all polynomials.  The Demazure atoms decompose the Schur functions in a natural way and their generating diagrams satisfy a Robinson-Schensted-Knuth-style algorithm~\cite{Mas08}.  Type A key polynomials~\cite{LasSch90, ReiShi95} are positive sums of Demazure atoms~\cite{Mas09}.  Summing the Demazure atoms over all weak compositions which collapse to a fixed composition when their zeros are removed produces a new collection of quasisymmetric functions, called the \emph{quasisymmetric Schur functions}, which we now formally define using fillings of composition diagrams ~\cite{HLMvW09}.

Let $\alpha$ be a composition of $n$.  If $T$ is a filling of the composition diagram $\alpha$ (written in English notation) satsifying the following properties, then $T$ is called a \emph{semi-standard reverse composition tableau}, abbreviated $\SSRCT$.

\begin{enumerate}
\item The entries in each row weakly decrease when read from left to right.
\item The entries in the leftmost column strictly increase when read from top to bottom.
\item (Triple Rule)  If $k>j$ and $T(i,k) \ge T(i+1,j)$ (for cells $(i,k)$ and $(i+1,j)$), then $(i+1,k)$ is a cell in $\alpha$ and $T(i+1,k) > T(i+1,j)$.  (Here, if there is no cell at coordinate $(i,j)$, set $T(i,j)=0$.)
\end{enumerate}

The set of all semi-standard reverse composition tableaux of shape $\alpha$ is denoted $\SSRCT(\alpha)$.  The \emph{weight} of a semi-standard reverse composition tableau $T$, denoted $X^T$, is the product over all $i$ of $x_i^{\#(i)}$, where $\#(i)$ is the number of times $i$ appears in $T$.

\begin{definition}
The \emph{quasisymmetric Schur function} $\mathcal{S}_{\alpha}$ is defined by $$\mathcal{S}_{\alpha}(X) = \sum_{T \in \SSRCT(\alpha)} X^T.$$
\end{definition}

Quasisymmetric Schur functions form a basis for $\QSym$ and are closely related to Schur funcions.  In fact, the quasisymmetric Schur functions, when summed over all rearrangements of a given partition, produce the Schur function indexed by this partition~\cite{HLMvW09}.  That is, $$s_{\lambda} = \sum_{\tilde{\alpha} = \lambda} \mathcal{S}_{\alpha},$$ where $\tilde{\alpha}$ is the partition obtained by arranging the parts of $\alpha$ into weakly decreasing order.

For example, the four semi-standard reverse composition tableaux of shape $(2,1)$ are $$\tableau{1 & 1 \\ 2}, \; \tableau{1 & 1 \\ 3}, \; \tableau{2 & 1 \\ 3}, \; {\rm and} \; \tableau{2 & 2 \\ 3},$$ producing the quasisymmetric Schur function $$\mathcal{S}_{21}(x_1,x_2,x_3) = x_1^2x_2+x_1^2x_3+x_1x_2x_3+x_2^2x_3.$$  The four semi-standard reverse composition tableaux of shape $(1,2)$ are $$\tableau{1 \\ 2 & 2}, \; \tableau{ 1 \\ 3 & 2}, \; \tableau{1 \\ 3 & 3}, \; {\rm and} \; \tableau{2 \\ 3 & 3},$$ producing the quasisymmetric Schur function $$\mathcal{S}_{12}(x_1,x_2,x_3) =x_1x_2^2 + x_1x_2x_3 + x_1x_3^2+x_2x_3^2.$$  Together these sum to $s_{21}(x_1,x_2,x_3)$; that is, $$s_{21} = \mathcal{S}_{21}+\mathcal{S}_{12}.$$

The quasisymmetric Schur functions expand positively in the fundamental basis for $\QSym$.  The set of all \emph{standard reverse composition tableaux} of shape $\alpha$, abbreviated $\SRCT(\alpha)$, is the subset of $\SSRCT(\alpha)$ consisting of the semi-standard reverse composition tableaux in which each of the positive integers in the set $\{1,2, \hdots, | \alpha | \}$ appears exactly once.  Each standard reverse composition tableau $T$ has a \emph{descent set} $\Des(T)$ defined by $$\Des(T) = \{ i \mid i+1 \; {\rm appears \; weakly \; right \; of} \; i \} \subseteq [n-1].$$

\begin{theorem}{\cite{HLMvW09}}
The quasisymmetric Schur function $\mathcal{S}_{\alpha}$ decomposes into the fundamental basis for quasisymmetric functions as follows: $$\mathcal{S}_{\alpha} = \sum_{T \in \SRCT(\alpha)} F_{\Des(T), |\alpha|}.$$
\end{theorem}

For example, the three standard reverse composition tableaux of shape $(2,1,3)$ are $$\tableau{3 & 1 \\ 4 \\ 6 & 5 & 2}, \qquad \tableau{2 & 1 \\ 4 \\ 6 & 5 & 3}, \qquad \text{and} \qquad \tableau{2 & 1 \\ 3 \\ 6 & 5 & 4}.$$  The descent sets are, respectively $\{1, 3, 4\}, \{2,4\},$ and $\{2,3\}$.  This implies that $$\mathcal{S}_{213} = F_{\{1,3,4\},6}+F_{\{2,4\},6} + F_{\{2,3\},6}.$$

Tewari and van Willigenburg~\cite{TewvWi15} introduce a collection of operators $\{\pi_i \}_{i=1}^{n-1}$ on standard reverse composition tableaux (which satisfy the same relations as the generators $\{ T_i \}_{i=1}^{n-1}$ described in Section~\ref{subsec:Hecke}) to produce an $H_n(0)$-action on standard reverse composition tableaux of size $n$.  

In particular, for $T \in \SRCT(\alpha)$ for some composition $\alpha \models n$ and $1 \le i \le n-1$, entries $i$ and $i+1$ are said to be \emph{attacking} if they are in the same column of $T$ or they are in adjacent columns of $T$ with $i+1$ appearing to the right of $i$ in a strictly lower row.  The operators $\pi_i$ for $1 \le i \le n-1$ are defined as follows, where $s_i(T)$ interchanges the positions of entries $i$ and $i+1$.

$\pi_i(T) = 
\begin{cases}
T &\textnormal{if} \; i \notin \Des(T), \\
0 & \textnormal{if} \; i \in \Des(T),  \; i \; {\rm and} \; i+1 \; \textnormal{are attacking, and} \\
s_i(T) & \textnormal{if} \;  i \in \Des(T),  \; i \; {\rm and} \; i+1 \; \textnormal{are non-attacking.}
\end{cases}$

Extend these operators to all of $\mathfrak{S}_n$ by setting $\pi_{\sigma} = \pi_{i_1} \pi_{i_2} \cdots \pi_{i_{\ell}}$ when $\sigma=s_{i_1} s_{i_2} \cdots s_{i_{\ell}}$ is any reduced word for $\sigma$.  Define a partial order $\preceq_{\alpha}$ on $\SRCT(\alpha)$ by setting $T_1 \preceq_{\alpha} T_2$ if and only if $\pi_{\sigma}(T_1) = T_2$ for some permutation $\sigma \in \mathfrak{S}_n$.  Extend $\preceq_{\alpha}$ to a total order $\preceq_{\alpha}^t$ arbitrarily and let $\mathcal{V}_{T_i}$ be the $\mathbb{C}$-linear span of all $T_j \in \SRCT(\alpha)$ such that $T_j \succeq_{\alpha}^t T_i$.

\begin{theorem}{\cite{TewvWi15}}
If $T_1 \in \SRCT(\alpha)$ is the minimal element under the total order $\preceq_{\alpha}^t$, then $\mathcal{V}_{T_1}:={\bf S_{\alpha}}$ is an $H_n(0)$-module whose quasisymmetric characteristic is the quasisymmetric Schur function $\mathcal{S}_{\alpha}$.
\end{theorem}

A \emph{simple} composition is a composition $\alpha=(\alpha_1, \alpha_2, \hdots, \alpha_{\ell})$ such that if $\alpha_i \ge \alpha_j \ge 2$ and $1 \le i < j \le \ell$, then there exists an integer $k$ satisfying $i < k < j$ such that $\alpha_k = \alpha_j-1$.

Tewari and van Willigenburg~\cite{TewvWi15} prove that ${\bf S_{\alpha}}$ is an indecomposable $H_n(0)$-module if and only if $\alpha$ is simple.  These results lead to the introduction of a new basis for quasisymmetric functions called the \emph{canonical quasisymmetric functions} $\{ \mathcal{C}_{\alpha} \}_{\alpha}$ and a branching rule for the ${\bf S_{\alpha}}$ which is analogous to the classical branching rule for Schur functions~\cite{Sag01}.

The product of a quasisymmetric Schur function and a Schur function expands into the quasisymmetric Schur function basis through a rule which refines the Littlewood-Richardson Rule~\cite{HLMvW11} but a formula for the coefficients appearing in the product of arbitrary quasisymmetric Schur functions is unknown.  See ~\cite{LMvW13} for a thorough introduction to quasisymmetric Schur functions and their closely related counterpart, the \emph{Young quasisymmetric Schur functions}.  The Young quasisymmetric Schur functions are obtained from the quasisymmetric Schur functions by a simple reversal of the indexing composition and the variables, but at times the Young quasisymmetric functions are easier to work with due to their compatibility with semi-standard Young tableaux (rather than reverse semi-standard Young tableaux). 

\section{Quasisymmetric analogues of symmetric function bases}{\label{sec:bases}}

Quasisymmetric functions play a major role in answering important questions about symmetric functions.  Analogues in $\QSym$ of classical bases for symmetric functions aid in this pursuit by providing a dictionary to translate between $\Sym$ and $\QSym$.  We have already discussed a quasisymmetric analogue of the monomial symmetric functions as well as two different quasisymmetric analogues of the Schur functions.  We now introduce another natural quasisymmetric analogue of the Schur function basis as well as a quasisymmetric analogue of the power sum basis.

\subsection{Dual immaculate quasisymmetric functions}

Berg, Bergeron, Saliola, Serrano, and Zabrocki~\cite{BBSSZ14} generalize Bernstein's creation operator construction of the Schur functions to obtain a basis for $\NSym$ called the \emph{immaculate basis} and denoted $\I_{\alpha}$.  The dual basis in $\QSym$, called the \emph{dual immaculate quasisymmetric functions}, can be generated by fillings of tableaux as follows.  

Let $F \colon \alpha \rightarrow \mathbb{Z}^+$ be a filling of a composition diagram $\alpha$ with positive integers such that the sequence of entries in each row (read from left to right) is weakly increasing and the sequence of entries in the leftmost column (read from bottom to top) is strictly increasing.  Then $F$ is said to be an \emph{immaculate tableau} of shape $\alpha$.  (Note that since we are using French notation our definition varies slightly from the definition in~\cite{BBSSZ14} but produces the same diagrams modulo a horizontal flip.)  The \emph{weight} of an immaculate tableau $U$, denoted $x^U$, is the product over all $i$ of $x_i^{\#(i)}$, where $\#(i)$ is the number of times $i$ appears in $U$.

\begin{definition}{\cite{BBSSZ14}}
Let $\alpha$ be a composition.  The \emph{dual immaculate quasisymmetric function} $\dI_{\alpha}$ is given by $$\dI_{\alpha}=\sum_U x^U,$$ where the sum is over all immaculate tableaux of shape $\alpha$.
\end{definition}

\begin{figure}
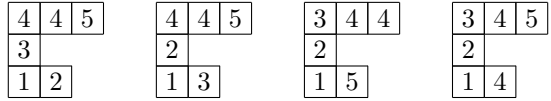

$$\tableau{4 & 4 & 5 \\ 3 \\ 1 & 2} \qquad \tableau{4 & 4 & 5 \\ 2  \\ 1 & 3} \qquad \tableau{3 & 4 & 4 \\ 2 \\ 1 & 5} \qquad \tableau{3 & 4 & 5 \\ 2 \\ 1 & 4}$$
\caption{The four immaculate tableaux of shape $(2,1,3)$ and weight $x_1x_2x_3x_4^2x_5$.}{\label{fig:dI}}
\end{figure}

For example, the coefficient of $x_1x_2x_3x_4^2x_5$ in $\dI_{2,1,3}$ is $4$ since there are four immaculate tableaux of shape $(2,1,3)$ and weight $x_1x_2x_3x_4^2x_5$.  (These immaculate tableaux are given in Figure~\ref{fig:dI}.)

The following theorem provides a formula for the expansion of the Schur functions into the dual immaculate quasisymmetric functions.

\begin{theorem}{\cite{BBSSZ14}}
Let $\lambda$ be a partition of length $k$.  Then $$s_{\lambda} = \sum_{\sigma} (-1)^{\sigma} \dI_{\lambda_{\sigma_1}+1 -\sigma_1, \lambda_{\sigma_2}+2-\sigma_2, \hdots , \lambda_{\sigma_k}+k-\sigma_k},$$ where $(-1)^{\sigma}$ is the sign of $\sigma$ and the sum is over all permutations $\sigma \in \symgrp_k$ such that $\lambda_{\sigma_i} + i - \sigma_i > 0$ for all $1 \le i \le k$.
\end{theorem}

For example, $$s_{321} = \dI_{321} - \dI_{141}.$$  Note that the coefficients are not always non-negative and further the compositions indexing the terms appearing in this expansion are not merely rearrangements of the partition $\lambda$ as is the case in the quasisymmetric Schur expansion of the Schurs.  However, the beauty of the connection to Schur functions is more readily apparent in the dual, since applying the forgetful map to an immaculate function produces the corresponding Schur function.  (That is, $\chi(\I_{\alpha}) = s_{\alpha}$.)

Grinberg recently proved Zabrocki's conjecture that the dual immaculate quasisymmetric functions can also be constructed using a variation on Bernstein's creation operators~\cite{Gri17}.  The dual immaculate quasisymmetric functions expand into positive sums of the monomial quasisymmetric functions, the fundamental quasisymmetric functions, and, recently shown in~\cite{AHM17}, the Young quasisymmetric Schur functions.  The latter expansion is not at all obvious given the very different methods used to generate these two bases, and therefore provides further justification that both of these families of functions are interesting and natural objects of study.

Like the quasisymmetric Schur functions, dual immaculate quasisymmetric functions correspond to characteristics of certain representations of the $0$-Hecke algebra~\cite{BBSSZ15}, but for the dual immaculate quasisymmetric functions these representations are indecomposable.  In particular, let $\mathcal{M}_{\alpha}$ be the vector space spanned by all words on the letters $\{1,2, \hdots , \ell(\alpha) \}$ such that the letter $j$ appears $\alpha_j$ times.  Define an action of the $0$-Hecke algebra on words by
$$\pi_i(w) = 
\begin{cases}
w & w_i \ge w_{i+1} \\
s_i(w) & w_i < w_{i+1},
\end{cases}$$
where $s_i(w)=w_1w_2 \cdots w_{i-1}  w_{i+1}w_i w_{i+2} \cdots w_n$.  Note that this is isomorphic to the induced representation $$\Ind_{H_{\alpha_1}(0) \otimes H_{\alpha_2}(0) \otimes \cdots \otimes H_{\alpha_{\ell(\alpha)}}(0)}^{H_n(0)} ( L_{\alpha_1} \otimes L_{\alpha_2} \otimes \cdots \otimes L_{\alpha_m}),$$ where $L_k$ is the one-dimensional representation indexed by the composition $(k)$.  A word $w$ in which the first instance of $j$ appears before the first instance of $j+1$ is called a $\mathcal{Y}$-word.  The $0$-Hecke action defined above cannot move a $j+1$ to the right of a $j$, so the subspace $\mathcal{N}_{\alpha}$ of $\mathcal{M}_{\alpha}$ spanned by all words which are not $\mathcal{Y}$-words is a submodule of $\mathcal{M}_{\alpha}$.

\begin{theorem}{\cite{BBSSZ15}}
The characteristic of $\mathcal{V}_{\alpha}:=\mathcal{M}_{\alpha} / \mathcal{N}_{\alpha}$ is the dual immaculate quasisymmetric function indexed by $\alpha$.  In other words, $\mathcal{F}([ \mathcal{V}_{\alpha}]) = \dI_{\alpha}$.
\end{theorem}

Bergeron, S\'anchez-Ortega, and Zabrocki found a Pieri rule (first conjectured in ~\cite{BBSSZ14} and proved in~\cite{BSZ16}) for the product of a fundamental quasisymmetric function and a dual immaculate quasisymmetric function, and much is known about the multiplication of the immaculate basis.  However, multiplication rules in full generality for the dual immaculate quasisymmetric functions are still largely unknown.

\subsection{Quasisymmetric analogues of the power sum basis}{\label{sec:powersum}}

The power sum symmetric functions (defined in Section~\ref{sec:background}) are eigenvectors for the omega involution $\omega$; that is, $\omega(p_{\lambda}) = \varepsilon_{\lambda} p _{\lambda}$, where $\varepsilon_{\lambda} = (-1)^{n- \ell(\lambda)}$~\cite{Sta99}.  Power sum symmetric functions are also helpful in computing characters of the symmetric group via the Murnaghan-Nakayama Rule~\cite{Mur37, Nak41}.

Malvenuto and Reutenauer~\cite{MalReu95}, through the Hopf algebraic dual, $\NSym$, of $\QSym$, introduce a quasisymmetric analogue of the power sum symmetric functions, also obtained independently by Derksen~\cite{Der09} using a similar process but with a computational error which leads to a different formula.  To understand their construction, we recall several facts about generating functions for symmetric and noncommutative symmetric functions.  The complete homogeneous symmetric functions, elementary symmetric functions, and power sum symmetric functions (in $n$ variables) can be defined through their generating functions $$H(t) = \sum_{d \ge 0} h_d t^d = \prod_{i=1}^n \frac{1}{1-x_it} \;  ,  \; E(t) = \sum_{k \ge 0} e_k t^k = \prod_{i=1}^n (1+x_it), \; {\rm and} \; P(t) = \sum_{k \ge 1} p_k \frac{t^k}{k}.$$  The relationship between these is given by Newton's formula:  $$-\frac{d}{dt} (E(-t)) = P(t)E(t),$$ which is equivalent to $$\frac{d}{dt}  ( H(t)) = H(t) P(t).$$  In their seminal work on noncommutative symmetric functions, Gelfand, Krob, Lascoux, Leclerc, Retakh, and Thibon~\cite{GKLLRT95} define a noncommutative analogue of the complete homogeneous symmetric functions (denoted ${\bf S}_k$) by describing their generating function (see Section~\ref{sec:Hopf}) and requiring they satisfy the multiplicative property.

They then utilize this approach to construct two analogues of the power sums in $\NSym$ by requiring that the generating functions satisfy the appropriate analogues of Newton's formula.  \emph{Noncommutative power sum symmetric functions of the first kind}, denoted ${\bf \Psi}_k$, are defined by $$\psi(t) = \sum_{k \ge 1} t^{k-1}{\bf \Psi}_k, \; \; \; \frac{d}{dt} \sigma(t) = \sigma(t) \psi(t), \; \; \; {\rm and} \; \; \; {\bf \Psi}_{\alpha} = {\bf \Psi}_{\alpha_1} {\bf \Psi}_{\alpha_2} \cdots {\bf \Psi}_{\alpha_{\ell}},$$ where $\sigma(t)$ is as defined in Equation~\ref{eqn:genfun}.  Similarly, \emph{noncommutative power sum symmetric functions of the second kind}, denoted ${\bf \Phi}_k$, are defined by $$\sigma(t) = \exp(\sum_{k \ge 1} t^k \frac{{\bf \Phi}_k}{k}) \; \; {\rm and} \; \; {\bf \Phi}_{\alpha} = {\bf \Phi}_{\alpha_1} {\bf \Phi}_{\alpha_2} \cdots {\bf \Phi}_{\alpha_{\ell}}.$$

Taking the Hopf algebraic duals of these noncommutative power sum bases produces two different quasisymmetric analogues of power sums.  We use ${\Psi}$ and ${ \Phi}$ as notation for these to emphasize their relationship with their noncommutative duals.  The dual of the noncommutative power sum basis of the first kind is defined~\cite{BDHMN17} by $${\Psi}_{\alpha} = z_{\alpha} \sum_{\beta \preceq \alpha} \frac{ M_{\beta}}{\pi(\alpha, \beta)},$$ where the ordering used is the refinement partial order (so that $\alpha \succeq \beta$ if $\alpha$ is coarser than $\beta$) and $\pi( \alpha, \beta)$ is given by the following process.  First define $\pi(\alpha)=\prod_{i=1}^{\ell(\alpha)} \sum_{j=1}^i \alpha_j$.  Then for $\alpha$ a refinement of $\beta$, set $\pi(\alpha,\beta)=\prod_{i=1}^{\ell(\beta)} \pi(\alpha^{(i)}),$ where $\alpha^{(i)}$ consists of the parts of $\alpha$ that combine to $\beta_i$. 

For example, ${ \Psi}_{312}=(1 \cdot 2 \cdot 3) ( \frac{1}{3 \cdot 1 \cdot 2} M_{312}+\frac{1}{3 \cdot 4 \cdot 2} M_{42} + \frac{1}{3 \cdot 1 \cdot 3} M_{33} + \frac{1}{3 \cdot 4 \cdot 6} M_6),$ which simplifies to $${\Psi}_{312}= M_{312} + \frac{1}{4} M_{42} + \frac{2}{3} M_{33} + \frac{1}{12} M_6.$$

Similarly, a formula for quasisymmetric power sums of the second kind is also given in terms of the monomial quasisymmetric functions.  $${ \Phi}_{\alpha} = \sum_{\alpha \succeq \beta} \frac{M_{\beta}}{f(\alpha, \beta)},$$ where the ordering used is again the refinement partial order, and the function $f(\alpha, \beta)$ is given by the following process.  Assume $\beta = ( \beta_1, \beta_2, \hdots, \beta_k)$.  Write $\alpha$ as a concatenation $\alpha^{(1)} \alpha^{(2)} \cdots \alpha^{(k)}$ of compositions $\alpha^{(i)}$ where $\alpha^{(i)} \models \beta_i$.  Then $f(\alpha, \beta) = \ell(\alpha^{(1)})! \cdots \ell(\alpha^{(k)})!$. 

For example, ${\Phi}_{312}=(1 \cdot 2 \cdot 3) ( \frac{1}{1 \cdot 1 \cdot 1} M_{312}+\frac{1}{2 \cdot 1} M_{42} + \frac{1}{1 \cdot 2} M_{33} + \frac{1}{6} M_6),$ which simplifies to $${\Phi}_{312}=6 M_{312} + 3 M_{42} + 3 M_{33} + M_6.$$

This formula differs from that of Malvenuto and Reutenauer~\cite{MalReu95} (who use the notation $P_{\alpha}$ instead of ${\Phi}_{\alpha}$) only by a constant.  This constant ensures that the $\Phi_{\alpha}$ refine the symmetric power sums so that $$p_{\lambda} = \sum_{\tilde{\alpha} = \lambda} \Phi_{\alpha},$$ which is not true for the $P_{\alpha}$.

The reader might wonder about the duals of the elementary and complete homogeneous symmetric functions.  In fact, the noncommutative complete homogeneous symmetric functions are dual to the monomial quasisymmetric functions, while the noncommutative elementary symmetric functions are dual to the ``forgotten'' basis for quasisymmetric functions, whose combinatorial structure is largely unknown.

Recall that the fundamental quasisymmetric functions satsify the following relationship to monomial quasisymmetric functions: $$F_{\alpha} = \sum_{\beta \preceq \alpha} M_{\beta},$$ where $\beta \preceq \alpha$ again means that $\beta$ is a refinement of $\alpha.$  Hoffman~\cite{Hof15} studied a variation on the fundamental basis called the \emph{essential quasisymmetric functions $E_{\alpha}$}, obtained by reversing the inequality in the above equation so that $$E_{\alpha} = \sum_{\beta \succeq \alpha} M_{\beta}.$$  Summing over all coarsenings of $\alpha$ is a natural thing to do because of what the antipode map does to the monomial quasisymmetric functions:  $$S(M_{\alpha}) = (-1)^{\ell(\alpha)} E_{\alpha}.$$  Multiplication in the essential basis follows the same rules (modulo a sign) as multiplication in the monomial basis.

\subsection{The shuffle algebra}

The shuffle algebra is a Hopf algebra (whose multiplicative structure is given by an operation called a \emph{shuffle}) which is in fact isomorphic as a graded Hopf algebra to $\QSym$ (over the rationals).  More details on the shuffle algebra and the closely related concept of Lyndon words can be found in Reutenauer~\cite{Reu93}, Lothaire~\cite{Lot97}, or Grinberg-Reiner~\cite{GriRei14}.  

Let $A$ be a totally ordered set, which we will call an \emph{alphabet}.  A \emph{word} of length $n$ is an ordered string $w=w_1w_2 \cdots w_n$ of elements of $A$.  Let $A^*$ be the set of all words on the alphabet $A$.  For the purposes of this section, we will take the alphabet to be the positive integers, as is done in~\cite{Haz01}.  When $A$ is taken to be the positive integers, the \emph{degree} ($|u|$) of a word $u$ in $A^*$ is the sum of its letters rather than the number of letters.  The shuffle, $u \shuffle v$, of two words $u=u_1 u_2 \cdots u_k$ and $v=v_1v_2 \cdots v_{\ell}$ in $A^*$ is the sum of all words  in $A^*$ of length $k+\ell$ formed from the letters of $u$ and the letters of $v$ such that for all $i$, $u_i$ appears before $u_{i+1}$ and $v_i$ appears before $v_{i+1}$. Multiplicities will occur if a letter appears in both $u$ and $v$.  If a letter appears more than once within one of the words $u$ or $v$, simply consider each occurrence as a distinct letter by applying a different subscript to each appearance of a given letter.  This product is associative and can therefore be extended to the shuffle product of a finite number of words.  (It is called a shuffle because it resembles the interleaving method used to shuffle a deck of cards.)  For example, $$23 \shuffle 12=2312 + 2132 + 2123 + 1223 +1223 +1232.$$

Shuffles in fact guide the multiplication of quasisymmetric power sums of both types.  Let $a_j$ equal the number of parts of size $j$ in $\alpha$ and $b_j$ equal the number of parts of size $j$ in $\beta$, and let $\alpha\cdot \beta$ denote their concatenation.  Define $C(\alpha,\beta)=\prod_j \binom{a_j+b_j}{a_j}$, so that $C(\alpha,\beta) = z_{\alpha\cdot\beta}/(z_\alpha z_\beta)$.

\begin{theorem}\label{thm:productpower}{\cite{BDHMN17}}
Let $\alpha$ and $\beta$ be compositions.  Then \[\Psi_\alpha \Psi_\beta = \frac{1}{C(\alpha,\beta)}\sum_{\gamma \in \alpha\shuffle\beta} \Psi_\gamma \; \; \qquad {\rm and} \; \; \qquad \Phi_\alpha\Phi_\beta = \frac{1}{C(\alpha,\beta)}\sum_{\gamma\in \alpha\shuffle\beta}\Phi_\gamma. \]
\end{theorem}

The \emph{shuffle algebra} $K \langle A \rangle$ (where $K$ is a commutative ring with unit) is the set of all linear combinations over $K$ of words on an alphabet $A$, endowed with this shuffle product.  There are a number distinct proofs that this algebra is isomorphic to $\QSym$, including those of Hazewinkel~\cite{Haz01} and Hazewinkel-Gubareni-Kirichenko~\cite{HGK10}.  Note that Theorem~\ref{thm:productpower} in fact implies that the shuffle algebra is isomorphic to $\QSym$.  To see this, distribute the $C(\alpha, \beta)$ in the first equation in Theorem~\ref{thm:productpower} so that: $$\frac{\Psi_{\alpha}}{z_{\alpha}} \frac{\Psi_{\beta}}{z_{\beta}} = \sum_{\gamma \in \alpha\shuffle\beta} \frac{\Psi_{\gamma}}{z_{\alpha\cdot\beta}}.$$  Then map from the shuffle algebra to $\QSym$ via the map $\alpha \mapsto \frac{\Psi_{\alpha}}{z_{\alpha}}$.  Extend this map linearly to an isomorphism between the shuffle algebra and $\QSym$. 

We now discuss the algebraic structure of $\QSym$.  First we shall see that $\QSym$ over the rationals is a polynomial algebra in the quasisymmetric power sums.  Then we describe Hazewinkel's polynomial generators for $\QSym$ over the integers.

Let $w \in A^*$ be a word on the alphabet $A$.  Then a \emph{proper suffix} of $w$ is a word $v \in A^*$ such that there exists a nonempty $u \in A^*$ such that $w=uv$.  A \emph{prefix} of $w$ is a word $u \in A^*$ such that $w=uv$.  Let $\le_{A}$ be a total ordering on $A^*$ defined as follows.  Let $u=u_1 u_2 \cdots u_k$ and $v=v_1 v_2 \cdots v_m$.  If $u$ is a prefix of $v$ then $u \le_{A} v$.  Otherwise let $j$ be the smallest positive integer such that $u_j \not= v_j$.  If $u_j>v_j$ then $u >_{A} v$.  Otherwise $u <_{A} v$.

\begin{definition}
A \emph{Lyndon word} is a nonempty word $w \in A^*$ such that every nonempty proper suffix $v$ of $w$ satisfies $w <_{A} v$.  Let $\mathcal{L}$ denote the set of all Lyndon words.
\end{definition}

For example, the words $1324, 1323,$ and $11213$ are Lyndon words while the words $4132, 3241, 2332,$ and $2233$ are not.  The shuffle algebra is freely generated over the rationals by the Lyndon words.  




\begin{theorem}{\cite{Rad79}}{\label{thm:rad}}
Every element of $K \langle A \rangle$ can be uniquely expressed as a polynomial in the Lyndon words.  In other words, the shuffle algebra $K \langle A \rangle$ is the polynomial algebra in the Lyndon words.
\end{theorem}
 
One can think of Theorem~\ref{thm:rad} (commonly known as \emph{Radford's Theorem}) as the statement that for any vector space basis whose elements are indexed by words in $A^*$ and whose multiplication is given by shuffles, each basis element can be written as a polynomial in basis elements indexed by Lyndon words.


For example, $w=321$ is not a Lyndon word, but $\Psi_{321}$ can be expressed as a polynomial in quasisymmetric power sums indexed by Lyndon words.  That is, $$\Psi_{321} = \Psi_1 \cdot \Psi_2 \cdot \Psi_3 - \Psi_{23} \cdot \Psi_1 - \Psi_3 \cdot \Psi_{12} + \Psi_{123}.$$  Radford's theorem implies that quasisymmetric power sums indexed by Lyndon words form an algebraically independent generating set  for $\QSym$ over the rationals.  (See~\cite{BDHMN17},~\cite{GriRei14},~\cite{MalReu95}, and~\cite{Reu93} for further details.)

We now describe Hazewinkel's polynomial generators for $\QSym$, which are indexed by a subset of Lyndon words.  The ring $\mathbb{Z}[x_1,x_2, \hdots ,]$ is endowed with a well-known $\lambda$-ring structure via $$\lambda_i(x_j) = 
\begin{cases}
x_j & i=1 \\
0 & i > 1,
\end{cases}$$ for $j=1,2, \hdots,$.  Define a total ordering on compositions called the \emph{wll-ordering} (weight, length, lexicographic) by:
\begin{enumerate}
\item If $|\alpha| > |\beta|$, then $\alpha >_{wll} \beta$.
\item If $|\alpha| = |\beta|$ and $\ell(\alpha) > \ell(\beta)$, then $\alpha >_{wll} \beta$.
\item If $|\alpha| = |\beta|, \ell(\alpha) = \ell(\beta),$ and $\alpha >_{lex} \beta$, then $\alpha >_{wll} \beta$.
\end{enumerate}

For example, $$523 >_{wll} 11213 >_{wll} 323 >_{wll} 143.$$

Hazewinkel proves~\cite{Haz10} that applying $\lambda_n$ to the monomial quasisymmetric function indexed by a Lyndon word $\alpha$ produces $$\lambda_n(M_{\alpha}) = M_{\alpha^{\star n}} + (smaller),$$ where $\alpha^{\star n}$ denotes concatenation of $\alpha$ with itself $n$ times and $(smaller)$ is a $\mathbb{Z}$-linear combination of monomial quasisymmetric functions which are wll-smaller than $\alpha^{\star n}$.  For example, $$\lambda_2(M_{(1,2)}) = M_{1212} + \text{ some subset of the set } V \cup W \cup Y,$$ where $$V=\{ \text{all words with weight }\le5\},$$ $$W=\{ \text{all words of weight }6\text{ and length } \le3\},$$ and $$Y=\{ 1122, 1113, 1131\}.$$

\begin{theorem}{\cite{Haz10}}{\label{thm:Haz}}
Let $eLYN$ be the set of all Lyndon words $u=u_1u_2 \cdots u_m$ such that $\gcd\{u_1, u_2, \hdots , u_m\}=1$.  Then the set $\{ \lambda_n(M_{u}) \}_{u \in eLYN}$ for all $n \in \mathbb{N}$ freely generates the ring of quasisymmetric functions over the integers.
\end{theorem}

Athough the monomial quasisymmetric functions are not multiplicative, Theorem~\ref{thm:Haz} provides a way to construct a multiplicative generating set.  Therefore, $\QSym$ is a polynomial algebra over the integers in the set $\{ \lambda_n(M_{\alpha}) \}_{\alpha \in eLYN}$.  Note that Theorem~\ref{thm:Haz} therefore implies the Ditters conjecture.


\section{Connections to symmetric functions and the polynomial ring}{\label{sec:sym}}

This section discusses several recent developments connecting quasisymmetric functions to important open problems within symmetric functions and the polynomial ring.  We focus our scope to three topics: chromatic quasisymmetric functions, transitions from $\QSym$ to $\Sym$, and liftings of $\QSym$ bases to the polynomial ring.  We will unfortunately not be able to address the Eulerian quasisymmetric functions~\cite{ShaWac10}, which are in fact symmetric despite their definition in terms of quasisymmetric functions.  See~\cite{ShaWac10} for a wonderful introduction to these fascinating objects of study, including the important definitions and theorems as well as the research avenues they introduce.  We also regretfully omit the recently developed theory of dual equivalence; see ~\cite{Ass15} for information about this new paradigm and how to use it.

\subsection{Chromatic quasisymmetric functions}

Let $G=(V,E)$ be a graph with vertices $V$ and edges $E$ and let $S$ be a subset of the positive integers $\mathbb{P}$.  A \emph{proper $S$-coloring} of $G$ is a function $\kappa \colon V \rightarrow S$ such that if two vertices $i$ and $j$ are adjacent (i.e. $\{i,j \} \in E$), then $i$ and $j$ are assigned different colors (i.e. $\kappa(i) \not= \kappa(j)$).  The \emph{chromatic number} $\chi(G)$ is the minimum number of colors (size of $S$) necessary to construct a proper $S$-coloring of $G$. 

It is natural to ask how many proper $\{1,2, \hdots , m \}$-colorings exist for a graph $G$; this number is denoted $\chi_G(m)$.  It is a non-negative integer when $m$ is a positive integer, and it is a polynomial called the \emph{chromatic polynomial} when $m$ is an indeterminant.  Stanley generalized this notion~\cite{Sta95} to construct a symmetric function generated from the set $\mathcal{C}(G)$ of all proper $\mathbb{P}$-colorings of $G$ as follows: $$X_G(x):= \sum_{\kappa \in \mathcal{C}(G)} x_{\kappa},$$ where $x=(x_1,x_2, \hdots )$ is a sequence of commuting variables and $x_{\kappa}=\prod_{v \in V} x_{\kappa(v)}$.  Notice that plugging in $x_i=1$ for all $i$ produces the chromatic polynomial $X_G(1^m) = \chi_G(m)$.  For example, the path $P_3$ on three vertices has chromatic symmetric function $$X_{P_3}(x_1,x_2,x_3) = x_1^2x_2+x_1^2x_3 + x_2^2x_3 + x_1x_2^2+x_1 x_3^2 + x_2x_3^2 + 6 x_1 x_2 x_3,$$ and $X_{P_3}(1,1,1)=12$, the number of proper colorings of $P_3$ with $3$ colors.

Recall that if a function has positive coefficients when expanded in a basis $B$, then it is said to be \emph{$B$-positive}.  For example, the elementary symmetric functions are Schur-positive since $$e_{\lambda} = \sum_{\mu} K_{\mu' \lambda} s_{\mu},$$ where $K_{\mu' \lambda}$ is the number of semi-standard Young tableaux of shape $\mu'$ and content $\lambda$.  One significant open question about chromatic symmetric functions relates to positivity in the elementary basis for symmetric functions.  If $P$ is a partially ordered set, then the \emph{incomparability graph} of $P$ is the graph $inc(P)$ whose vertices are the elements of $P$ and whose edges are the pairs of vertices which are incomparable in $P$.  A poset is called \emph{$(r+s)$-free} if no induced subposet is isomorphic to the direct sum of a chain (totally ordered set) with $r$ elements and a chain with $s$ elements.

\begin{conjecture}{(Stanley-Stembridge Conjecture~\cite{Sta95,StaSte93})}
If $G=inc(P)$ for some $(3+1)$-free poset $P$, then $X_G(x)$ is $e$-positive.
\end{conjecture}

\begin{figure}
\begin{center}
\xymatrix{
a \ar@{-}[rd] & & c  & d &   & & & a \ar@{-}[r] \ar@{-}[rd] & c \\
 & b \ar@{-}[ru] & &  & & & & b  \ar@{-}[r] & d \ar@{-}[u] \\
& &  P & & & & & & inc(P)
}
\end{center}
\caption{A poset $P$ and its incomparability graph}
{\label{fig:graphs}}
\end{figure}
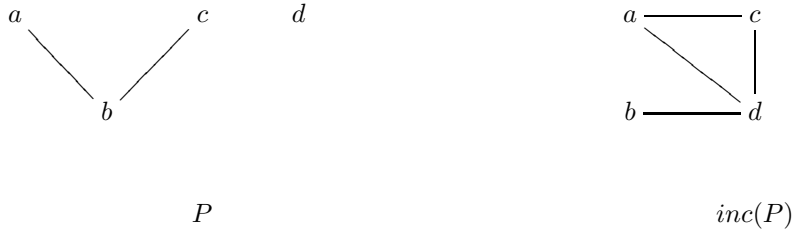

For example, the poset $P$ in Figure~\ref{fig:graphs} is $(3+1)$-free.  The $e$-expansion for the chromatic symmetric function corresponding to its incomparability graph is $X_{inc(P)} = 4e_{31}+8e_{4}$.

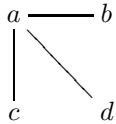
\begin{figure}
\begin{center}
\hspace*{.05in} 
\xymatrix{
a \ar@{-}[r] \ar@{-}[rd] \ar@{-}[d] & b \\
c & d
}
\end{center}
\caption{The star graph $S_3$ on four vertices}
{\label{fig:star}}
\end{figure}

The incomparability graph for a $(3+1)$-free poset is an example of a \emph{claw-free} graph.  A claw-free graph is a graph which does not contain the star graph $S_{3}$ (depicted in Figure~\ref{fig:star}) as a subgraph.  However, not all claw-free graphs are $e$-positive; Dahlberg, Foley and van Willigenburg~\cite{DFvW17} provide a family of claw-free graphs which are not $e$-positive.

Gasharov~\cite{Gas96} proved that the incomparability graph of a $(3+1)$-free poset is Schur-positive.  Since the elementary symmetric functions are Schur-positive, Schur-positivity would follow immediately from $e$-positivity.  Guay-Paquet~\cite{Gua13} proved that the chromatic symmetric function of a $(3+1)$-free poset is a convex combination of chromatic symmetric functions of posets which are both $(3+1)$-free and $(2+2)$-free.  This reduces the $e$-positivity conjecture to a subclass of posets with more structure than posets which are $(3+1)$-free.  A \emph{natural unit interval order} is a poset $P$ on the set $[n]:=\{1,2, \hdots, n\}$ obtained from a certain type of intervals on the real line as follows.  Let $\{ [a_1, a_1+1], [a_2, a_2+1], \hdots , [a_n, a_n+1] \}$ be a collection of closed intervals of length one such that $a_i < a_{i+1}$ for $1 \le i \le n-1$.  Set $i <_P j$ if $a_i+1 < a_j$.  The resulting partially ordered set will always be $(3+1)$-free and $(2+2)$-free, and in fact every poset that is both $(3+1)$-free and $(2+2)$-free is a unit interval order~\cite{ScoSup58}.  

Shareshian and Wachs recently proposed a new approach to the Stanley-Stembridge $e$-positivity conjecture in the form of a refinement of Stanley's chromatic symmetric functions.  This refinement behaves nicely with respect to unit interval orders. 

\begin{definition}{\cite{ShaWac16}}
Let $G=(V,E)$ be a graph whose vertex set $V$ is a finite subset of $\mathbb{P}$.  The \emph{chromatic quasisymmetric function} of $G$ is $$X_G(x,t) = \sum_{\kappa \in \mathcal{C}(G)} t^{\asc(\kappa)} x^{\kappa},$$ where $$\asc(\kappa) = | \{ \{i,j \} \in E \mid i<j \; {\rm and} \; \kappa(i) < \kappa(j) \} |.$$
\end{definition}

\begin{figure}
\begin{center}
\hspace*{.05in} 
\xymatrix{
1 \ar@{-}[r] \ar@{-}[rd] \ar@{-}[d] & 2 \\
4 & 3
}
\end{center}
\caption{A labeling of the star graph $S_3$ on four vertices}
{\label{fig:chromstar}}
\end{figure}
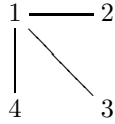

Notice that the chromatic quasisymmetric function $X_G(x,t)$ depends not only on the isomorphism class of the graph $G$ but also on the labeling of the vertices of $G$.  Let $G$ be the star graph on four vertices labeled as shown in Figure~\ref{fig:chromstar}.  Then $$X_{G}(x,t) = M_{31}+M_{121} + M_{211} + M_{1111} + t(M_{121}+2M_{211}+3M_{1111})+$$ $$ + t^2(M_{121}+2M_{211}+M_{112}+5M_{1111}) + t^3 (M_{121}+2M_{112}+M_{211}+M_{13}+6M_{1111}) +$$ $$+ t^4(M_{121}+2M_{112}+5M_{1111})+t^5(M_{121}+M_{112}+3M_{1111})+t^6(M_{1111}).$$  The chromatic quasisymmetric function reduces to the chromatic symmetric function by setting $t=1$; that is $X_G(x,1)=X_G(x)$.

Let $\omega$ be the involution map on $\QSym$ which sends $F_{S}$ to $F_{[n-1] \setminus S}$.  The image of $X_G(x,t)$ under $\omega$ has a natural positive expansion into the fundamental basis for quasisymmetric functions~\cite{ShaWac16}.  Shareshian and Wachs~\cite{ShaWac12} further conjecture that when $G$ is the incomparability graph of a natural unit interval order, this image corresponds to a sum of Frobenius characteristics associated to certain Hessenberg varieties.  This conjecture was proved by Brosnan and Chow~\cite{BroCho15} and, through a different approach, by Guay-Paquet~\cite{Gua16}, providing an alternate proof of Schur positivity.  


\begin{theorem}{\cite{ShaWac16}}
If $G$ is the incomparability graph of a natural unit interval order, then $X_G(x,t)$ is symmetric in the $x$-variables.
\end{theorem}

Not every graph whose chromatic quasisymmetric function is symmetric is an incomparability graph of a natural unit interval order.  One interesting open question is to classify which graphs admit a symmetric chromatic quasisymmetric function.

Several extensions of chromatic quasisymmetric functions have recently emerged, demonstrating the many different areas this research impacts.  Ellzey extends this paradigm to directed graphs~\cite{Ell17}.  Haglund and Wilson express the integral form Macdonald polynomials as weighted sums of chromatic quasisymmetric functions~\cite{HagWil17}.  Clearman, Hyatt, Shelton, and Skandera interpret the chromatic quasisymmetric functions in terms of Hecke algebra traces~\cite{CHSS16}, while Alexandersson and Panova connect the chromatic quasisymmetric functions to LLT polynomials~\cite{AlePan17}.

\subsection{Quasisymmetric expansions of symmetric functions}

As quasisymmetric functions become more ubiquitous, many natural expansions of symmetric functions into quasisymmetric functions (particularly into the fundamental quasisymmetric functions) are appearing.  It is natural to try to use this structure to answer classical questions about symmetric functions such as Schur positivity.  Egge, Loehr, and Warrington~\cite{ELW10} recently introduced a method to convert the quasisymmetric expansion of a symmetric function into the Schur function expansion, providing a new approach to questions of Schur positivity.  

We need several definitions in order to describe the ``modified inverse Kostka matrix'' and some interesting applications of this paradigm.  A \emph{rim-hook} is a set of contiguous cells in a partition diagram such that each diagonal contains at most one cell.  A \emph{special rim-hook tableau} is a decomposition of a partition diagram into rim-hooks such that each rim-hook contains at least one cell in the leftmost column of the diagram.  E\u{g}ecio\u{g}lu and Remmel~\cite{EgeRem90} use special rim-hook tableaux in their formula for the \emph{inverse Kostka matrix}, which is the transition matrix from the monomial basis for symmetric functions to the Schur functions.

The \emph{sign} of a special rim-hook is $(-1)^{r-1}$, where $r$ is the number of rows spanned by the rim-hook.  The sign of a special rim-hook tableau is the product of the signs of its rim-hooks.  A special rim-hook tableau is said to be \emph{flat} if each rim-hook contains exactly one cell in the leftmost column of the partition diagram.

\begin{theorem}{\cite{ELW10}}{\label{thm:ELW}}
Let $F$ be a field, and let $f$ be a symmetric function given by its expansion into the fundamental quasisymmetric functions so that $$f=\sum_{\alpha \models n} y_{\alpha} F_{\alpha}.$$  Then the coefficients $x_{\lambda}$ in the Schur function expansion $\displaystyle{f=\sum_{\lambda \vdash n} x_{\lambda} s_{\lambda}}$ are given by $$x_{\lambda} = \sum_{\alpha \models n} y_{\alpha} K^*_n(\alpha, \lambda),$$ where $K^*_n(\alpha,\lambda)$ is the sum of the signs of all flat special rim-hook tableaux of partition shape $\lambda$ and content $\alpha$.
\end{theorem}

Theorem~\ref{thm:ELW} provides a potential alternative approach to proving that Macdonald polynomials expand positively into the Schur functions.  In particular, recall that Theorem~\ref{thm:HHL} describes a formula for expanding Macdonald polynomials into the fundamental quasisymmetric functions.  Combining this formula with Theorem~\ref{thm:ELW} implies that the coefficient of $s_{\lambda}$ in the Schur function expansion of $\tilde{H}_{\mu}$ is given by 
$$\sum_{\alpha \vdash n} K^*_n(\alpha, \lambda) \left( \sum_{\beta \in \{\mathfrak{S}_n \mid \Des(\beta^{-1})=\alpha\} } q^{\mathrm{\inv}(\beta,\mu)} t^{\mathrm{\maj}(\beta,\mu)} \right).$$

The following example is similar to that appearing in~\cite{ELW10}.  If $\mu=(3,1)$ and $\lambda=(2,2)$, then there exists a flat special rim-hook tableau for $\alpha=(2,2)$ and a flat special rim-hook tableau $\alpha=(1,3)$.  These correspond to $$K^*_4((2,2),(2,2))=+1 \textnormal{ and } K^*_4((1,3),(2,2))=-1$$ respectively.  The permutations $w$ whose inverse descent sets $\Des(w^{-1})$ are $\{2\}$ are $3412,3142,3124,1324,$ and $1342$.  Computing the $\inv$ and $\maj$ for the fillings of $(3,1)$ with these permutations as reading words produces $2q^2+qt+t+q$.  Similarly, the permutations whose inverse descent sets are $\{ 1 \}$ are $2341, 2314,$ and $2134$.  Their $\inv$ and $\maj$ (for fillings of $(3,1)$) produce $-q^2-q-t$.  Putting this together, the coefficient of $s_{22}$ in $H_{31}$ is $$2q^2+qt+t+q - (q^2+q+t) = q^2+qt.$$

Notice that negative terms do appear in the $K^*_n(\alpha, \lambda)$.  This means that in order to apply this technique to the Schur positivity of Macdonald polynomials problem, one must find involutions to cancel out the negative terms.

A further application of this transition matrix from the fundamental quasisymmetric functions to Schur functions is to the Foulkes Plethysm Conjecture~\cite{Fou50}, which states that $s_n[s_m]-s_m[s_n]$ (where the brackets denote a certain type of substitution called \emph{plethysm}) is Schur positive.  Loehr and Warrington~\cite{LoeWar12} provide a formula for the expansion of $s_{\mu}[s_{\nu}]$ into fundamental quasisymmetric functions using a novel interpretation of the ``reading word'' of a matrix.  The modified inverse Kostka matrix could then be used to determine the Schur function expansion of $s_n[s_m]-s_m[s_n]$, again with the caveat that involutions are needed to cancel out the negative terms.  

Garsia and Remmel recently found a further extension of the Egge, Loehr, Warrington result.  They proved that each fundamental appearing in the fundamental expansion of a symmetric function can be replaced by the Schur function indexed by the same composition.  Since every such Schur function is either $0$ or $\pm s_{\lambda}$ for some partition $\lambda$, this expansion can be simplified to a signed sum of Schur functions indexed by partitions.

\begin{theorem}{\cite{GarRem18}}{\label{thm:GarRem}}
Let $f$ be a symmetric function which is homogeneous of degree $n$ and expands into the fundamental basis for quasisymmetric functions as follows:
$$f = \sum_{\alpha \models n} a_{\alpha} F_{\alpha}.$$  Then $$f = \sum_{\alpha \models n} a_{\alpha} s_{\alpha}.$$
\end{theorem}

Theorem~\ref{thm:GarRem} already has a number of important consequences.  Garsia and Remmel used this approach to formulate a conjecture regarding the modified Hall-Littlewood polynomials.  Leven applied this method to prove an extension of the Shuffle Conjecture for the cases $m=2$ and $n=2$~\cite{Lev14}.  Qiu and Remmel considered the cases of this ``Rational Shuffle Conjecture'' where $m$ or $n$ equals $3$~\cite{QiuRem17}.  

\subsection{Slide polynomials and the quasi-key basis}

Schubert polynomials are an important class of polynomials, first introduced by Lascoux and Sch$\ddot{{\rm u}}$tzenberger~\cite{LasSch82} to provide a new method for computing intersection numbers in the cohomology ring of the complete flag variety.  Several different combinatorial formulas for Schubert polynomials have been discovered since their original introduction as divided difference operators, including but not limited to reduced decompositions~\cite{BJS93},~\cite{FomSta94} and $RC$-graphs~\cite{BerBil93}.  Despite the numerous ways to construct Schubert polynomials, it remains an open problem to provide a combinatorial formula for the expansion of a product of Schubert polynomials into the Schubert basis.

Assaf and Searles~\cite{AssSea17} further the study of Schubert polynomials with the introduction of two new families of polynomials, both of which positively refine the Schubert polynomials.  These new families, called the \emph{monomial slide polynomials} and the \emph{fundamental slide polynomials}, exhibit positive structure constants (meaning the coefficients appearing in their products are always positive), whereas the key polynomials (another family of polynomials refining the Schubert polynomials~\cite{Dem74a,LasSch90,ReiShi95}) have signed structure constants.  Although the slide polynomials have many interesting applications (to Schubert polynomials and other objects of study in algebraic combinatorics), this article focuses on their connections to quasisymmetric functions.  

Remove the zeros from a weak composition $\gamma$ to obtain a (strong) composition called the \emph{flattening} of $\gamma$, denoted $\flat(\gamma)$.  The \emph{monomial slide polynomial} $\mathfrak{M}_{\gamma}$ is then defined by $$\mathfrak{M}_{\gamma}(x_1, x_2, \hdots , x_n) = \sum_{\substack{\delta \ge \gamma \\ \flat(\delta)=\flat(\gamma)}} x_1^{\delta_1}x_2^{\delta_2} \cdots x_n^{\delta_n},$$ where $\delta \ge \gamma$ if $\delta$ \emph{dominates} $\gamma$; that is $\delta_1+\delta_2 + \cdots + \delta_i \ge \gamma_1 + \gamma_2 + \cdots + \gamma_i$ for all $1 \le i \le n$.  The related \emph{fundamental slide polynomial} $\mathfrak{F}_{\gamma}$ is defined by $$\mathfrak{F}_{\gamma} = \sum_{\substack{\delta \ge \gamma \\ \flat(\delta) \; {\rm refines} \; \flat(\gamma)}} x_1^{\delta_1}x_2^{\delta_2} \cdots x_n^{\delta_n}.$$  For example, \begin{align*} \mathfrak{F}_{1032}(x_1,x_2,x_3,x_4) &= x_1x_3^3x_4^2+x_1x_2^3x_4^2+x_1x_2^3x_3^2+x_1x_2x_3^2x_4^2+x_1x_2^2x_3x_4^2+x_1x_2^3x_3x_4 \\ &=\mathfrak{M}_{1032}+\mathfrak{M}_{1122} + \mathfrak{M}_{1212} + \mathfrak{M}_{1311}.\end{align*}  Each of these families $\{ \mathfrak{M}_{\gamma} \}_{\gamma}$ and $\{ \mathfrak{F}_{\gamma} \}_{\gamma}$ (indexed by weak compositions of $k$) of polynomials is a $\mathbb{Z}$-basis for polynomials of degree $k$ in $n$ variables.  

Assaf and Searles~\cite{AssSea16} also introduce a related basis, called the \emph{quasi-key polynomials} $\mathfrak{Q}_{\gamma}$, for the polynomial ring which is analogous to the key polynomials.  These polynomials are positive sums of fundamental slide polynomials and in fact stabilize to the quasisymmetric Schur functions as zeros are prepended to their indexing compositions.  (Prepending $m$ zeros to the composition $\gamma$ is denoted by $0^m \times \gamma$.)

\begin{theorem}{\cite{AssSea16}}
For any weak composition $\gamma$, we have $$\lim_{m \rightarrow \infty} \mathfrak{Q}_{0^m \times \gamma} = \mathcal{S}_{\flat(\gamma)}.$$
\end{theorem}

Each Schubert polynomial can be written as a positive sum of fundamental slide polynomials using a new object called a \emph{quasi-Yamanouchi pipe dream}.  While this definition takes us too far from our current topic, we do take the time to describe a closely related construction involving the fundamental expansion of Schur functions.

Recall that the Schur functions decompose into a positive sum of the fundamental quasisymmetric functions (see Equation~\ref{Schur:fund}); this formula can be computed by finding the descent sets of all standard Young tableaux of a given shape.  However, when the number of variables is less than the number of descents, the corresponding fundamental equals $0$.  Assaf and Searles~\cite{AssSea17} introduce a class of semi-standard Young tableaux, called \emph{quasi-Yamanouchi tableaux}, which dictate precisely which fundamentals appear in the decomposition with non-zero coefficient when the variables are restricted.  

\begin{definition}{\cite{AssSea17}}
A semi-standard Young tableau is said to be \emph{quasi-Yamanouchi} if for all $i>1$, the leftmost occurrence of $i$ lies weakly left of some appearance of $i-1$.  Let $QYT_n(\lambda)$ denote the set of quasi-Yamanouchi tableaux of shape $\lambda$ whose entries are in $[n]$.
\end{definition}

The \emph{weight} of a quasi-Yamanouchi tableau $T$ is given by $wt(T) = \prod_i x_i^{m_i}$, where $m_i$ is the number of times the entry $i$ appears in $T$. 

\begin{theorem}{\cite{AssSea17}}
The Schur polynomial $s_{\lambda}(x_1, \hdots , x_n)$ is given by $$s_{\lambda}(x_1, \hdots , x_n) = \sum_{T \in QYT_n(\lambda)} F_{wt(T) }(x_1, \hdots , x_n).$$
\end{theorem}

For example, the three quasi-Yamanouchi tableaux of shape $\lambda=(4,2)$ and entries in $\{1,2\}$ are $$\tableau{2 & 2 \\ 1 & 1 & 1 & 1},  \; \; \tableau{2 & 2 \\ 1 & 1 & 1 & 2}, \;  \; \tableau{2 & 2 \\ 1 & 1 & 2 & 2 },$$ and therefore the Schur expansion into fundamentals is $$s_{42}(x_1,x_2) = F_{4,2}(x_1,x_2)+F_{3,3}(x_1,x_2)+F_{2,4}(x_1,x_2).$$

Note that all the terms appearing on the right hand side are nonzero, and there is no need to calculate the descent sets for all nine standard Young tableaux of shape $(4,2)$.  This is important in the study of Schubert polynomials because although certain classes of Schubert polynomials are equal to Schur functions, the number of variables appearing varies based on the indexing permutation.  For example, the Schubert polynomial indexed by the permutation $213$ (written in one-line notation) is equal to the Schur function $s_1(x_1) = x_1$ while the Schubert polynomial indexed by the permutation $132$ (written in one-line notation) is $s_1(x_1,x_2) = x_1 + x_2$.  A thorough understanding of precisely the nonzero terms appearing in the quasisymmetric expansion is therefore crucial to the quest of proving a combinatorial formula for Schubert multiplication.

This connection to Schubert multiplication (a long-standing open problem in algebraic combinatorics) exemplifies the utility of quasisymmetric functions.  Quasisymmetric functions appear in a number of other important problems which have helped to shape the study of algebraic combinatorics including Schur positivity of Macdonald polynomials, the Foulkes plethysm conjecture, and the Stanley-Stembridge conjecture.  We hope the reader comes away from this article with a deeper appreciation for the beauty and utility of quasisymmetric functions and a desire to further explore this exciting and far-reaching avenue of research.

\section{Acknowledgements}\label{sec:ack}

I am very grateful to H\'{e}l\`{e}ne Barcelo, Gizem Karaali, and Rosa Orellana for inviting me to produce this chapter.  I would also like to thank Ed Allen, Susanna Fishel, Josh Hallam, Jim Haglund, and John Shareshian for helpful feedback along the way.  Finally, I greatly appreciate the insightful comments from a diligent anonymous referee.

\bibliographystyle{plain}
\bibliography{Qsymbib}
\label{sec:biblio}

\end{document}